\documentclass[12pt]{amsart}

\usepackage[colorlinks=true, pdfstartview=FitV, linkcolor=blue, citecolor=blue, urlcolor=blue]{hyperref}

\usepackage{amssymb,amsmath, amscd}
\usepackage{times, verbatim}
\usepackage{graphicx}
\input xy
\xyoption{all}

\usepackage{mathrsfs}
\usepackage[mathcal]{euscript}
\usepackage{dsfont}
\usepackage[T1]{fontenc}

\usepackage[small,nohug,heads=littlevee]{diagrams}
\usepackage[normalem]{ulem}

\usepackage{enumerate}
\usepackage{epsfig}
\usepackage{MnSymbol}
\RequirePackage{color}
\usepackage{mathrsfs}
\usepackage{tikz}
\usepackage{tikz-cd}

\usetikzlibrary{matrix,arrows}

\usepackage{amsfonts}
\usepackage{amscd}
\usepackage{mathrsfs}
\usepackage{amsmath}
\usepackage{enumerate}
\usepackage{epsfig}
\RequirePackage{color}
\usepackage{soul}

\usepackage{pgf}
\usepackage{tikz}
\usepackage{tikz-cd}

\usetikzlibrary{matrix,arrows}

\definecolor{purple}{rgb}{0.59, 0.44, 0.84}
\newcommand{\ft}[1]{{\color{blue}  #1}}

\definecolor{turk}{rgb}{0.2, .60, .50}

\newtheorem{theorem}{Theorem}[section]

\newtheorem{proposition}[theorem]{Proposition}

\newtheorem{example}[theorem]{Example}

\def \Tr {\mr{Tr}}
\def \f12{\frac12}
\def \Frob{\mr{Frob}}

\def\CCe#1#2{\binom {#1}{#2}}

\def\C{\mathbb C}

\def\P{\mathbb P}

\newcommand{\ZZ}{\mathbf{Z}}
\newcommand{\FF}{\mathbf{F}}
\newcommand{\RR}{\mathbf{R}}

\newcommand{\QQ}{\mathbf{Q}}

\newcommand{\PP}{\mathbb{P}}
\newcommand{\CC}{\mathbf{C}}

\newcommand{\qlbar}{\bar{\QQ}_{\ell}  }

\newcommand{\mc}[1]{\mathcal{#1}}
\newcommand{\mbf}[1]{\mathbf{#1}}
\newcommand{\mfr}[1]{\mathfrak{#1}}
\newcommand{\mr}[1]{\mathrm{#1}}

\newcommand{\slr}{\mr{SL}_2( \RR)}

\newcommand{\slz}{\mr{SL}_2( \ZZ)}

\newcommand{\G}{\Gamma}

\def\({\left(}
\def\){\right)}

\def\l{\lambda}
\def\eps{\varepsilon}
\def \ol{\overline}

\newcommand*\HYPERskip{&}
\catcode`,\active
\newcommand*\pFq{
	\begingroup
	\catcode`\,\active
	\def ,{\HYPERskip}%
	\doHyper
}
\catcode`\,12
\def\doHyper#1#2#3#4#5{%
	\, _{#1}F_{#2}\left[\begin{matrix}#3 \smallskip \\  #4\end{matrix} \; ; \; #5\right]%
	\endgroup
}

\catcode`,\active

\catcode`\,12

\catcode`,\active

\catcode`\,12

\newcommand*\HYPERpp{&}
\catcode`,\active
\newcommand*\pPPq{
	\begingroup
	\catcode`\,\active
	\def ,{\HYPERpp}%
	\doHyperFpp
}
\catcode`\,12
\def\doHyperFpp#1#2#3#4#5{%
	\, _{#1}{\mathbb P}_{#2}\left[\begin{matrix}#3 \smallskip \\  #4\end{matrix} \; ; \; #5\right]%
	\endgroup
}

\catcode`,\active

\catcode`\,12

\catcode`,\active

\catcode`\,12

\catcode`,\active

\catcode`\,12

\begin{document}
	\title{QM abelian varieties, hypergeometric character sums and modular forms}
	
	\author{Jerome William Hoffman}
	
	\address{Department of Mathematics \\
		Louisiana State University \\
		Baton Rouge, Louisiana 70803}

	\author{Fang-Ting Tu}
	
	\address{Department of Mathematics \\
		Louisiana State University \\
		Baton Rouge, Louisiana 70803}

	\email{hoffman@math.lsu.edu}
	
	\email{tu@math.lsu.edu}
	
	
	
	\begin{abstract}
		This is a report on recent work, with Wen-Ching Winnie Li and Ling Long. In that work explicit 
		formulas are given, involving  hypergeometric character sums, 
		for the traces of Hecke operators $T_p$  acting spaces
		of cusp forms $S_k(\Gamma)$ of weight $k$ for certain 
		arithmetically defined Fuchsian subgroups
		$\Gamma\subset \mr{SL}_2 (\RR)$. In particular we consider the groups attached 
		to the quaternion division algebra $B_6$ over $\QQ$ of discriminant 6.

	\end{abstract}

	\maketitle
	
	\section{Introduction.}
	\label{S:intro}
	The forthcoming joint work with Wen-Ching Winnie Li and Ling Long \cite{HLLT}  extends some investigations of  Wen-Ching Winnie Li, Tong Liu, and  Ling Long on
	4-dimensional Galois representations with quaternion structures, \cite{ALLL}; see also 
	\cite{LLT2}. Our main results are  explicit formulas in the shape (see Section \ref{S:main}; terms to be explained in more detail)
	\begin{equation}
		\label{E:main}
		-\mr{Tr}(T_p \mid S_k(\Gamma)) =
		\text{an expression in hypergeometric character sums.}
	\end{equation}
	Here $\mr{Tr}(T_p \mid S_k(\Gamma))$ is the trace of the Hecke operator $T_p$  acting on the space
	of weight-$k$ cusp forms $S_k(\Gamma)$ on the arithmetic Fuchsian group $\G$. 
	In fact, the expression will be a sum of terms for each $\lambda \in  X_{\Gamma} (\FF _p)  $ of  the local Frobenius traces
	$\mr{Tr}(\mr{Frob} _{\lambda})$ on a constructible
	$\bar{\QQ}_{\ell}$-sheaf $ V^k (\Gamma) $
	for the \'etale topology on \ft{the} curve $X_{\Gamma}$. 
	The theory of  hypergeometric character sums over finite 
	fields was  developed 
	largely by Greene \cite{Greene},  Katz \cite{Katz}, Beukers-Cohen-Mellit \cite{BCM}, 
	and  Fuselier-Long-Ramakrishna-Swisher-Tu \cite{Win3X}. The groups $\Gamma$ are certain 
	arithmetic triangle groups.

	The inspiration for this work can be found in the papers of Ahlgren, Frechette, Fuselier, Lennon, Ono, 
	Papanikolas
	\cite{Ahlgren01, AO00,FOP2, Fuselier, lennon2, lennon1}, wherein
	formulas in terms of hypergeometric character sums for traces of Hecke operators
	are given, but not in terms of traces of Frobenius in \'etale cohomology. 
	Rather, their approach combines 
	(1) counting formulas for elliptic 
	curves over finite fields due to Schoof, \cite{schoof}, and (2) the Selberg trace formula. In fact, formulas for the 
	trace of Hecke operators $T_p$ on  $S_{2k}(\mr{SL}_2 (\ZZ))$ had already been given by Ihara, \cite{Ihara67} by
	similar methods, but Ihara did not use hypergeometric character sums. 
	
	We give a different approach, which is more geometric, and applies not only to the
	cases treated by these authors, but also to groups $\Gamma$ arising from the units in the quaternion algebra
	$B_6$ over $\QQ$ of discriminant 6. We will explain the geometric viewpoints on modular forms and 
	hypergeometric functions.

	\section{Background.}
	\label{S:back}
	This work weaves several major mathematical threads. On the one hand, 
	this is part of the general Langlands program, namely that part concerned with expressing the 
	$L$-functions of {\it motives} in terms of $L$-functions of {\it automorphic forms}. In the case at hand, 
	the automorphic forms are classical modular forms on the algebraic group $GL (2)_{\QQ}$ or its inner 
	twist $B_{6, \QQ}$ given by the quaternion algebra $B_6$. The automorphic $L$-functions are 
	attached to {\it local systems on modular curves}. When  $\Gamma\subset \mr{SL}_2 (\ZZ)$ is a 
	congruence subgroup, these are the modular curves classifying families of elliptic curves with level structure. 
	When $\Gamma$ comes from a quaternion algebra, these are Shimura curves. Thus, the motives involved 
	belong to the general theory of {\it Shimura varieties}, \cite{Milne05}. 
	
	Another major theme is that of hypergeometric functions. The term {\it hypergeometric} 
	includes not only the classical $\phantom{}_2 F _1$ functions, and their generalizations
	$\phantom{}_{p} F _{q} $, Appell-Lauricella functions, etc. There are (at least) two general extensions 
	of the theory of hypergeometric functions: (1) the $A$-hypergeometric, or GKZ (Gelfand-Kapranov-Zelevinsky)
	systems, \cite{GKZ90}; (2) the Gabber-Loeser-Sabbah hypergeometric systems, based on earlier works of M. Sato and Ore, 
	\cite{LS312, LS91, LS315}. See the paper of Dwork/Loeser: \cite{DworkLoeser93}.

	Generally speaking, a hypergeometric function is one that 
	
	\begin{itemize}
		\item[1.] has power-series expansions in special form: $\Gamma$-series;
		\item[2.] satisfies a (regular) holonomic system of differential equations;
		\item[3.] has Euler integral expressions;
		\item[4.] is attached to a motivic sheaf.
	\end{itemize}
	
	What the last item means is that not only are there complex-analytic functions defined by them, but also $\ell$-adic
	and $p$-adic versions. The $\ell$-adic versions give rise to hypergeometric character sums, which are central to this 
	paper. For the GKZ systems see \cite{LF}; for the GLS systems see \cite{GL1}. The $p$-adic versions give $p$-adic analytic functions, 
	first investigated by Dwork. For the GKZ systems, see \cite{FWZ}.

	Note that there are irregular differential equations of hypergeometric type, the confluent hypergeometrics. 
	The Hodge-deRham realizations are related to irregular Hodge theory and the periods belong to exponential modules, 
	see \cite{DworkLoeser94}. The character  sums involve 
	additive as well as multiplicative characters of finite fields, and hence their $\ell$-adic sheaves have wild ramification at infinity. 
	For a recent example see \cite{FSY22}.

	\section{Modular forms and modular curves}
	\label{S:geom}
	
	\subsection{Modular forms}
	\label{SS:geom1}
	(Reference: \cite{Shimura-book}).
	A Fuchsian subgroup of the first kind is a discrete subgroup  $\Gamma \subset \slr$
	such that the invariant volume $\Gamma \backslash \slr$ is finite. This acts
	on the upper half-plane $\mathfrak{H}$ by
	\[
	\gamma . z = \frac{a z + b}{c z + d}, \quad \gamma = 
	\begin{pmatrix} a & b \\ c &d \end{pmatrix}.
	\]
	The upper half-plane is conformally equivalent to the open unit disk $\mathfrak{D}$, 
	via the map $z \mapsto w = (z-i)/(z+i)$.
	This reflects the isomorphism $\mr{SU}(1, 1) \cong \slr$, where
	\[
	\mr{SU}(1, 1) = \left \{ 
	\begin{pmatrix}  a & b  \\ \bar{b} & \bar{a}\end{pmatrix} \mid a, b \in \CC, \ \ |a|^2 - |b|^2 = 1.
	\right \}, \quad
	\begin{pmatrix}  a & b  \\ \bar{b} & \bar{a}\end{pmatrix} .w = \frac{a w + b}{\bar{b}w+ \bar{a}}.
	\]
	There are isomorphisms $\slr/\mr{SO}(2)\cong \mathfrak{H\cong }\mathfrak{D}$
	where 
	\[
	\mr{SO}(2) = \left \{\begin{pmatrix}    \cos(\theta)  & \sin(\theta)\\ -\sin(\theta) & \cos(\theta)\end{pmatrix}
	\right \}\subset \slr \text{\ \ is a maximal compact subgroup.}
	\]
	For a Fuchsian subgroup $\Gamma$ of first kind,
	let $ X_{\Gamma} ^{an} = \Gamma \backslash \mathfrak{H}$ or 
	$\Gamma \backslash \mathfrak{H}^*$ according to whether $\Gamma$ cocompact or not.
	Here $\mathfrak{H}^*$ is the union of $\mathfrak{H}$ and the cusps 
	of $\Gamma$, if any. We have a finite set
	$S = S_{c} \cup S_{e} \subset X_{\Gamma} $ of points which are cusps or elliptic points, 
	that is,  their preimages in $\mathfrak{H}^*$ are cusps or elliptic points of $\Gamma$, respectively. Let
	\[
	X_{\Gamma} ^{\circ} = X_{\Gamma} - (S_{c} \cup S_{e} ), \quad
	\quad
	Y_{\Gamma} = X_{\Gamma} - S_{c}  =  \Gamma \backslash \mathfrak{H}, 
	\]
	the complement of the set of cusps and elliptic points, and the complement of the 
	set of cusps, respectively. Note that 
	$p: \mathfrak{H} -\{  \text{elliptic points of \ }  \Gamma  \}\to  X_{\Gamma} ^{\circ} $
	is a covering space in the sense of topology. In particular, if $\Gamma$ has no elliptic points,
	then $p:  \mathfrak{H}  \to   Y_{\Gamma} $ is the universal covering of $ Y_{\Gamma}$.
	The action of $\Gamma$ on  $\mathfrak{H}^*$  is via the quotient 
	$\bar{\Gamma}$ in  $\mr{PSL}_2( \RR)$. For the groups in this paper,
	$\bar{\Gamma} = \Gamma /\Gamma \cap \{ \pm I\}$.  Therefore, if 
	$\Gamma$ torsion-free, then 
	$p:  \mathfrak{H}  \to   Y_{\Gamma} $ is the universal covering of $ Y_{\Gamma}$, and
	the fundamental group at any base-point $x$
	\[
	\pi _1 (Y_{\Gamma}, x) \cong \Gamma, 
	\]
	an isomorphism unique up to inner automorphism. In the quaternion cases 
	this is $\pi _1 (X_{\Gamma}, x) \cong \Gamma$ for $\Gamma$ without torsion.
	In general, there is an epimorphism 
	\[
	\pi _1 (X^{\circ}_{\Gamma}, x)\to \bar{\Gamma } \cong \mr{Aut} (( \mathfrak{H} -\{  \text{elliptic points of \ }  \Gamma  \})/ X _{\Gamma}^{\circ}   ).
	\]

	If $\Gamma \subset \slr$ is a Fuchsian subgroup of the first kind, and $k\ge 0$ is an integer, we recall that a modular form 
	of weight $k$ is a holomorphic function $f(z)$ for $z \in \mathfrak{H}$ such that 
	\[
	f(\gamma .z) = j(\gamma, z)^k f(z) \text {\ for all\ } \gamma \in\Gamma, \quad
	j\left (\begin{pmatrix} a & b \\ c &d \end{pmatrix} , z \right ) = cz+d
	\]
	which satisfies a growth condition at the cusps. 
	For each cusp $c$ of $\Gamma$ there is a parameter $q_c$ such that a modular form has 
	an expansion at $c$: $f = a_0 + a_1 q_c + a_2 q_c^2 +...$.  A cusp form is one that vanishes at every cusp
	i.e., $a_0 = 0$ at every cusp of $\Gamma$. The space of cusp forms $S_k (\Gamma)$ is finite 
	dimensional over $\CC$.  We assume the reader is familiar with the standard 
	subgroups $\Gamma (N), \Gamma _0(N) \subset \slz$ as well as the Hecke operators 
	$T(p)$ (the latter only defined in the case of arithmetically defined Fuchsian subgroups).

	\subsection{Modular curves}
	\label{SS:geom2}
	In our work, we consider those $\Gamma$ which are arithmetically defined subgroups. 
	There are two cases: 
	\begin{itemize}
		\item[1.] Elliptic modular case. $\Gamma$ is commensurable with  $\mr{SL}_2( \ZZ)$. Then the cusps of $\Gamma$ is the set 
		of rational numbers and the point $\infty$. 
		\item[2.]Quaternion case.  $\Gamma$ is commensurable to the set  $O_B ^1$ of norm 1 elements in a maximal order $O_B$ of an indefinite quaternion algebra $B$
		over $\QQ$. That is, we choose once and for all an embedding $B \subset M_2 (\RR)$, which induces an embedding 
		$\theta : O_B ^1 \subset \mr{SL}_2(\RR)$. In this case, the set of cusps is empty, so  $\mathfrak{H}^* = \mathfrak{H}$. 
	\end{itemize} 
	
	In both these cases, Shumura's theory of canonical models shows that these Riemann surfaces are the 
	$\CC$-points of an algebraic curve, denoted $X_{\Gamma}$ defined over a number field $k_{\Gamma}$.	
	We will use the notation $X_{\Gamma}$ for the corresponding Riemann surface, if there is no ambiguity. If necessary, we 
	use  $X ^{an}_{\Gamma}$ to denote the analytic space. 
	
	These curves are (in general coarse) moduli spaces. If $\Gamma$ has no torsion, then there are universal 
	families of abelian varieties $f: A_{\Gamma} \to Y_{\Gamma}$ with additional structures. In the elliptic modular cases these are families 
	of elliptic curves $f: E_{\Gamma} \to Y_{\Gamma}$. In the quaternion cases, 
	$Y_{\Gamma}= X_{\Gamma}$, these are families of 2-dimensional abelian varieties
	with an action of the quaternion algebra, that is, with a homomorphism
	\[
	\theta : B \to \mr{End}(A)_{\QQ}
	\] 
	and additional structure. 
	
	\begin{example}
		\label{E:M4}
		(See \cite{DR,Shioda73}). Let
		\[
		M_4 = \mr{Spec}\, \ZZ[i, 1/2, \sigma, (\sigma(\sigma ^4 -1))^{-1}].
		\]
		This is the moduli scheme for $\Gamma (4)\subset \slz$. 
		The universal 
		elliptic curve for this is 
		\[
		E_{\sigma}: y^2 =  x(x-1)(x-\lambda), \quad \lambda = (\sigma + \sigma ^{-1})^2/4.
		\]
		This curve is isomorphic with the Jacobi quartic
		\[
		C_{\sigma}: y^2 = (1-\sigma ^2 x^2)(1-x^2 /\sigma^2).
		\]
		via the change of variables
		\[
		X = \frac{ \sigma ^2+1}{2 \sigma ^2}\cdot\frac{x-\sigma}{x-1/\sigma}, \ \ 
		Y = \frac{\sigma ^4-1}{4 \sigma ^3}\cdot\frac{y}{(x-1/\sigma)^2}.
		\]
	\end{example}
	
	\begin{example}
		\label{M3} Let 
		\[
		M_3 = \mr{Spec}\, \ZZ[\zeta _3, 1/3, \mu, (\mu ^3 -1)^{-1}].
		\]
		This is the moduli scheme for $\Gamma (3)\subset \slz$. The universal 
		elliptic curve for this is 
		\[
		E_{\mu}: x^3 +y^3+z^3 -3 \mu xyz =0.
		\]
		
	\end{example}
	The quaternion case with $D=6$ will be discussed in detail in section \ref{S:D=6}.
	
	\subsection{Triangle groups}
	\label{SS:geom3}
	
	In this work, we consider those $\Gamma$ that are (or are closely related to) triangle groups.
	We recall the definition: Let $a, b, c \in \ZZ _{\ge 2}\cup \{ \infty \}$, with $a \le b \le c$. Define 
	\[
	\chi (a, b, c) =  \frac{1}{a}+ \frac{1}{b}+ \frac{1}{c}-1
	\]
	and refer to $(a, b, c)$ as {\it spherical},  {\it Euclidean},  {\it hyperbolic}, depending on whether 
	$\chi(a, b, c)$ is $>0$, $=0$ or $<0$ respectively. In each of these three cases we associate a geometry
	\[
	H = 
	\begin{cases}
		\mr{the\  sphere \ } \mathbb{P} ^1 (\CC), \text{\ \ \ \ \ \ \ \ \ \ \ \ \ \ \ \ if \ }\chi (a, b, c) > 0;\\
		\mr{the\  plane \ } \CC, \text{\ \ \ \ \ \ \ \ \ \ \ \ \ \ \ \ \ \ \ \ \ \ \ \ \  if \ }\chi (a, b, c) = 0;\\
		\mr{the\ upper\ half-plane \ } \mathfrak{H}, \text{\ \ \ \ if \ }\chi (a, b, c) < 0.
	\end{cases}
	\]
	Spherical triples: $(2, 3, 3)$, $(2, 3, 4)$, $(2, 3, 5)$, $(2, 2, c)$, for $c \ge 2$. Euclidean triples:
	$(2, 2, \infty)$, $(2, 3, 6)$, $(2, 4, 4)$, $(3, 3, 3)$. The rest are hyperbolic. Among them, 
	there are finitely many arithmetic triangle groups, which were enumerated by Takeuchi, see \cite{Takeuchi-triangle}, 
	\cite{Takeuchi-classify}. To each triple we have 
	the associated triangle group, defined as 
	\[
	\Delta (a, b, c) = \langle  \delta_a, \delta_b, \delta _c \mid \delta _a ^a =   \delta _b ^b =  \delta _c ^c = 
	\delta _a \delta _b \delta _c  = 1\rangle.
	\]
	We have representation of $\Delta (a, b, c) $ via isometries of the associated geometry $H$ as follows. 
	Let $T$ be a geodesic triangle in $H$ with angles $\pi /a, \pi/b, \pi/c$. Let $\tau _a$, $\tau _b$, $\tau _c$
	be reflections in the three sides of $T$. The group generated by $\tau _a, \tau _b, \tau _c$ is a discrete group 
	with fundamental domain $T$. The subgroup of orientation-preserving isometries is generated by 
	\[
	\delta _a  =  \tau _c \tau _b, \quad \delta _b  =  \tau _a \tau _c,  \quad \delta _c  =  \tau _b \tau _a 
	\]
	satisfies the relations of $\Delta (a, b, c)$ with $\delta_p$  a counterclockwise rotation at the vertex with angle $2\pi i/p$. The quotient space $X(\Delta(a,b,c)):=\Delta(a,b,c)\backslash H$ is a Riemann orbifold of genus zero.   The fundamental domain $D_{\Delta}$ is then a union of $T$ and
	a reflected image of $T$.

	\begin{figure}[h]
		\begin{center}
			\includegraphics[scale=.5]{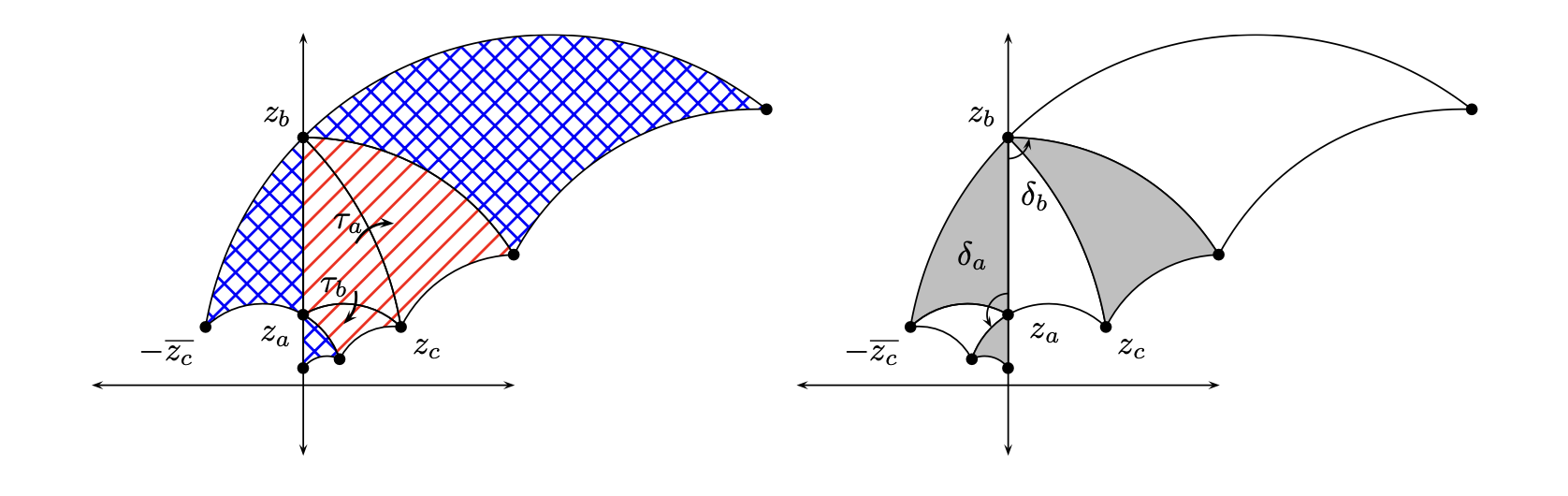}
		\end{center}
		\caption{A normalized triangle and its images under the reflections 
			$\tau_a $  and $\tau_b$, and under the rotations $\delta _a  = \tau_c \tau _b$
			and $\delta _b = \tau_a \tau _c$. Image from Klug, Musty, Schiavone and Voight, \cite{KMSV14}. }
	\end{figure}

	The relevance of triangle groups for hypergeometric equations arises from  the fundamental theorem of H. A. Schwarz
	\cite{Schwarz73}. 
	Given two independent solutions $y_1$ and $y_2$ to a hypergeometric differential equation, 
	the ratio $\eta  = y_1 /y_2$ maps the complex upper half-plane conformally onto a region $T$ 
	in the extended complex plane which is bounded by circular arcs. For the function
	\[
	\pFq{2}{1}{a,  b}{,c}{z} 
	\]
	the angles are $p \pi$, $q \pi $, $r \pi$ where $p = |1-c|$, $q = |c-a-b|$, $r = |a-b|$. 
	Schwarz observed that the monodromy of the differential equation gave rise to isometries of the associated 
	geometry $H$. In some cases, one obtains by analytic continuation a tessellation of $H$. 
	For this to extend to a tesselation we must have that $p, p, r$ are all in the shape $1/n$ where $n \in \ZZ_{\ge 2} \cup \{ \infty \}$. 
	We can denote the angles 
	of this triangle by  $\pi /a, \pi/b, \pi/c$ and classify the triples $(a, b, c)$ as above.
	He was especially interested 
	in the spherical case, where the monodromy group is necessarily finite, and he managed to completely classify these. 
	
	Given a triangle group $\G=(e_0,e_1,e_\infty)$,  following Theorem  9  of  \cite{Yang-Schwarzian} by Y. Yang we introduce
	the following  hypergeometric parameters: $$a=\frac12(1-\frac 1{e_1}-\frac1{e_0}-\frac1{e_\infty}),\quad b=\frac12(1-\frac 1{e_1}-\frac1{e_0}+\frac1{e_\infty}), \quad c=1-\frac 1{e_0}, $$
	and
	$$\tilde a=\frac12(1-\frac 1{e_1}+\frac1{e_0}-\frac1{e_\infty}),\quad \tilde b=\frac12(1-\frac 1{e_1}+\frac1{e_0}+\frac1{e_\infty}), \quad \tilde c=1+\frac 1{e_0}. $$ Using these Yang wrote down an explicit basis for $S_k(\G)$ in terms of $\pFq21{a&b}{&c}{t}$ and $\pFq21{\tilde a&\tilde b}{&\tilde c}{t}$.

	\begin{figure}[h]
		\begin{center}
			\includegraphics[scale=.3]{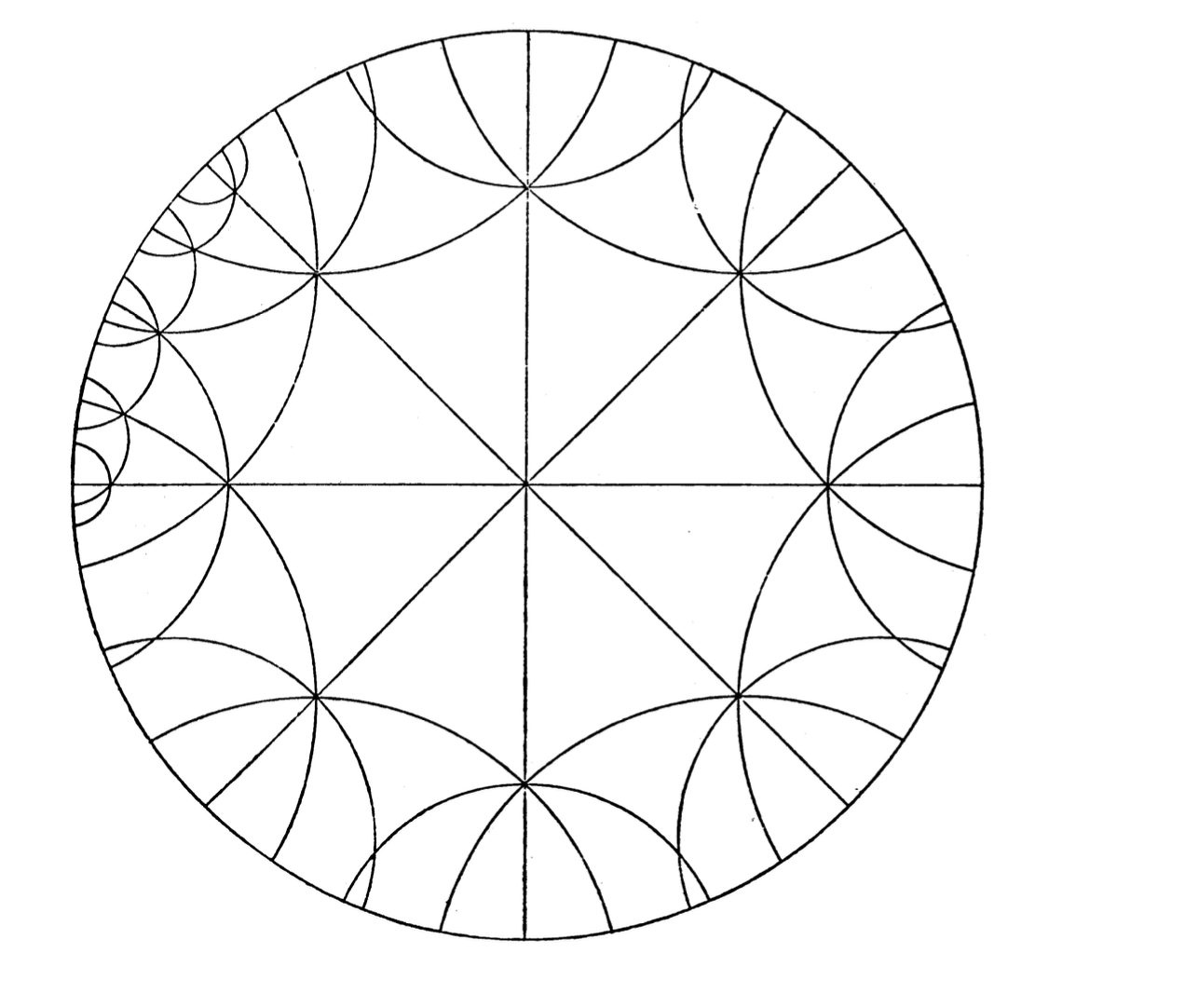}
		\end{center}
		\caption{ Tesselation of the hyperbolic disk 
			belonging to $(4, 4, 4)$ that appears in Gauss's collected works volume 8. }
	\end{figure}
	There are now extensive resources on the web for exploring hyperbolic tesselations.
	This had several major later developments. 
	On the one hand, Poincar\'e  extended this analysis to conformal maps given by second-order linear differential equations with $n$ regular singular points. The conformal image is an $n$-gon bounded by circular arcs. By studying 
	the transformations of these under analytic continuation, Poincar\'e made the link with non-Euclidean geometry and
	initiated the modern theories of hyperbolic geometry and automorphic forms. For hypergeometric differential equations of 
	several variables, including those studied by Picard, and Appell and Lauricella, see the works of Deligne and Mostow, 
	\cite{DM86}.

	\section{Local systems.}
	\label{S:Local}
	\subsection{Complex analytic.}
	\label{SS:Hyplocal1}
	
	Riemann introduced 
	the idea of monodromy into the study of analytic differential equations. 
	Given  a representation of the fundamental group
	\[
	\rho : \pi _1 (\PP^1 (\CC) - \{ 0, 1, \infty\}, x) \to \mr{GL}_2(\CC)
	\]
	there is a unique second order rational differential equation with regular singular points 
	$z = 0, 1, \infty$ with the property that if $f_1(z), f_2(z)$ is a basis of holomorphic 
	solutions at $x$ then analytic continuation around a loop $\gamma \in \pi _1 (\PP^1 (\C) - \{ 0, 1, \infty\}, x)$
	yields the linear transformation
	\[
	\begin{pmatrix} f_1\\ f_2 \end{pmatrix} \to \rho (\gamma) \begin{pmatrix} f_1\\ f_2 \end{pmatrix}.
	\]
	Since 
	\[
	\pi _1 (\PP^1 (\CC) - \{ 0, 1, \infty\}, x) = \langle \gamma _0, \gamma _1, \gamma _{\infty}\mid
	\gamma _0 \gamma _1 \gamma _{\infty} = 1\rangle , 
	\]
	to give the monodromy representation is equivalent to giving two-by-two matrices
	$T_0 , T_1, T_{\infty}  $ such that 
	$T_0  T_1 T_{\infty}= 1 $. These are well-defined up to 
	simultaneous conjugation by an element of  $\mr{GL}_2(\CC)$. As Katz observes, Riemann proved 
	a stronger result. Namely, it suffices to give the Jordan canonical forms of 
	$T_0 , T_1, T_{\infty}  $ to reconstruct the hypergeometric
	differential equation. Actually Riemann only considered the case when these were semisimple, so 
	equivalent to diagonal matrices, with eigenvalues $\exp (2\pi i \alpha )$ and  $\exp (2\pi i \alpha ' )$; 
	he called the  $\alpha, \alpha '$ the exponents at the singular point. They are well-defined up to permutation and
	adding $\ZZ$. This stronger property, that a differential equation is determined by the Jordan forms
	of the monodromies at the singular points, is what Katz (\cite{Katz96}) calls {\it rigidity}. Rigidity plays an important 
	role in our work.
	
	If $X$ is a connected nonsingular algebraic variety over $\CC$ and 
	\[
	\rho : \pi _1 (X^{an},  x) \to \mr{GL}_n(\CC)
	\]
	is a representation, we get a local system {\sf  V} on the analytic space $X ^{an}$. By the Riemann-Hilbert
	correspondence, this is the solution sheaf to a differential equation, unique up to isomorphism,
	\[
	\nabla : \mc{V} \to \Omega^ 1 _{X/\CC} \otimes _{\mc{O}_X} \mc{V}
	\]
	with regular singular points at infinity (\cite{DelDE}, \cite{Katz76over}). If $X$ is a nonsingular algebraic variety over $\CC$ we let $X^{an} = 
	X(\CC)$ be the set of complex points with the classical topology. Assume
	that $ X$ is quasi-projective and connected. Recall the following dictionary:
	The  following categories are equivalent: 
	\begin{itemize}
		\item[1.] Local systems of finite-dimensional $\CC$-vector spaces {\sf V }on 
		$ X^{an}$.
		\item[2.] Representations $\rho : \pi _1 (X^{an},  x) \to \mr{GL}(V)$
		on finite-dimensional $\CC$-vector spaces $V$.
		\item[3.] Holomorphic integrable connections 
		\[
		\nabla : \mc{V}^{an} \to \Omega^ 1 _{X^{an}/\CC} \otimes _{\mc{O}_{X^{an}}} \mc{V}^{an}
		\] 
		where $\mc{V}^{an}$ is  a locally free $\mc{O}_X^{an}$-module of finite rank.  
		\item[4.] Integrable algebraic connections 
		\[
		\nabla : \mc{V}\to \Omega^ 1 _{X/\CC} \otimes _{\mc{O}_{X}} \mc{V}
		\] 
		where $\mc{V}$ is  a locally free $\mc{O}_X$-module of finite rank, and which have regular 
		singular points ``at infinity''.   
	\end{itemize}
	
	\subsection{$\ell$-adic.}
	\label{SS:Local2}
	
	A reference: \cite{DelWeil1}, \cite{DelWeil2}. Fix  a prime number $\ell$.  In this section: scheme = 
	a separated noetherian scheme on which $\ell$ is invertible. We are interested in constructible $\bar{\QQ}_{\ell}$-sheaves on $X$, in particular, 
	those that are lisse. In this section: the \'etale topology is understood. An $\ell$-adic representation of a profinite group $\pi$ on a  $\bar{\QQ}_{\ell}$-vector 
	space $V$  is a homomorphism $$\sigma : \pi \to \mr{GL}(V)$$ such that 
	there is a finite subextension $E/\QQ _{\ell}$ and an $E$-structure $V_E$ on $V$
	such that $\sigma $ factorizes in a continuous homomorphism $\pi \to \mr{GL}(V_E)$.
	Recall that a geometric point $\bar{x}$ of a scheme $X$ is a morphism of the 
	spectrum of an algebraically closed field denoted $k(\bar{x})$ to $X$. It is localized 
	in $x \in X$ if its image is $x$. If $X$ is connected and pointed by a geometric point $\bar{x}$, the functor
	\[
	\mc{F} \mapsto \text{the \ } \pi_1 (X, \bar{x}) -\text{\ module\ } \mc{F}_{\bar{x}}
	\]
	is an equivalence between the categories of
	\begin{itemize}
		\item[1.] lisse  $\bar{\QQ}_{\ell}$-sheaves on $X$;
		\item[2.] $\ell$-adic continuous representations of $\pi _1(X, \bar{x})$.
	\end{itemize}
	Here   $\pi _1(X, \bar{x})$ is Grothendieck's fundamental group. Especially if $X = \mr{Spec}(k)$ is a field, the category of lisse  
	$\bar{\QQ}_{\ell}$-sheaves on $X$ is 
	equivalent to the category of $\ell$-adic representations of $\mr{Gal}(\bar{k}/k)$. 
	
	If one is working in the category of schemes $X$ of finite type over a perfect  field $k$ with algebraic closure
	$\bar{k}$, and $x \in X(k)$ is a $k$-rational point, 
	by convention we let $\bar{x}$ be the geometric point of $X$ which is the composite
	$\mr{Spec}(\bar{k}) \to \mr{Spec}(k) \to X $. 
	
	More generally, let $S$ be an irreducible scheme of finite type over a field  $k$. Then 
	a lisse sheaf $\mathcal{V}$ on $S$  is equivalent to a representation
	\[
	\rho : \pi_1 (S, \bar{\eta}) \to GL(\mathcal{V}  _{\bar{\eta}}) \sim GL(n, \bar{\QQ} _{\ell} )
	\] 
	where  $\bar{\eta}$ is a geometric generic point. That is,  
	$\bar{\eta} = \mathrm{Spec} (\overline{k (S) })$ where $\eta$ is the generic point
	of $S$, and  $k(S) = \mc{O}_{S, \eta}$ is the function field of $S$.

	Recall that there is a canonical exact sequence
	\[
	\begin{CD}
		0 @>>> \pi_1 (S _{\bar{k}}, \bar{\eta}) @>>>  \pi_1 (S , \bar{\eta})  @>>> \mathrm{Gal}(\bar{k}/k) @>>> 0.
	\end{CD}
	\]
	Moreover if $k$ has characteristic $0$ and we choose an embedding 
	$k \hookrightarrow \CC$,  $\pi_1 (S _{\bar{k}}, \bar{\eta}) $  is isomorphic to the profinite completion of the Poincar\'e fundamental 
	group, $\pi _1 (S^{an}, x) $ 
	so this includes the geometric monodromy of $\mathcal{F}$.
	
	Now let $k = \FF_q$ be a be a finite field of characteristic $p> 0$. 
	If $X_0$ is a scheme of finite 
	type, we  
	define $X= X_0 \otimes _{\FF _q}\bar{\FF}_q$ and recall that there is a morphism $F: X \to X$
	which sends a point with coordinates $x$ to the point with coordinates $x^q$. 
	The set of closed points $|X|$ is canonically identified with 
	$X_0(\bar{\FF}_q)$, and the Frobenius fixed points $|X|^F = X({\FF _q})$.
	The set of closed points $|X_0|$ is isomorphic to the set of orbits 
	$|X|_F$ of $F$; for $x_0 \in |X_0|$ the size of the corresponding orbit $Z$
	is $\mr{deg}(x_0)$, which is also the degree of the extension field
	$k(x)/\FF_q$. 
	If $\mc{F}_0$ is a constructible $\bar{\QQ}_{\ell}$-sheaf on $X_0$ for the \'etale topology, 
	we let $(X, \mc{F})$ be deduced from $(X_0, \mc{F}_0)$ by extension to the algebraic
	closure. If $\bar{x}\in |X|$ is a geometric point localized in $x \in |X_0|$, 
	we define
	\[
	\mr{Frob}_x := \mr{Frob}_{q^{\mr{deg}(x)}}\in \mr{Gal}(k(\bar{x})/k(x)) = \mr{Gal}(\bar{\FF} _q/\FF _q) .
	\]
	The stalk $\mc{F} _{\bar{x}}$ is a  $\bar{\QQ}_{\ell}$-vector space of finite dimension on which
	$\mr{Gal}(k(\bar{x})/k(x))$ acts and therefore 
	\[
	\det \left ( 1- t.   \mr{Frob}_x \mid     \mc{F} _{\bar{x}} \right )\in \bar{\QQ}_{\ell}[t]
	\] 
	is defined. It is independent of the geometric point $\bar{x}$ localized in $x$, so we can 
	simply denote it by $\det \left ( 1- t.   \mr{Frob}_x \mid     \mc{F}_0  \right )$. The coefficient 
	of $-t$ in this polynomial is the trace $\mr{Tr} (\mr{Frob}_x \mid \mc{F}_0)$. The Grothendieck-Lefschetz
	trace formula reads:
	\begin{equation}
		\label{E:GL}
		\sum _{i = 0}^{2d} (-1)^i \mr{Tr} (\mr{Frob}_q \mid H^i _c (X, \mc{F}   )) = 
		\sum _{x \in X(\FF_q)}\mr{Tr} (\mr{Frob}_x \mid \mc{F}_0), 
	\end{equation}
	where $d$ is the dimension of $X$ and $H^i _c$ is cohomology with proper support.  
	We define the $L$-function
	\[
	Z(X_0, \mc{F}_0, t):= \prod _{x \in |X_0|} \det \left ( 1- t ^{\mr{deg}(x)}   \mr{Frob}_x \mid     \mc{F}_0  \right )^{-1}.
	\]
	The trace formula is equivalent to the factorization
	\[
	Z(X_0, \mc{F}_0, t):= \prod _{i = 0}^{2d} \det \left ( 1- t \mr{Frob}_q \mid H^i _c (X, \mc{F})\right )^{i+1}.
	\]

	\subsection{Families of varieties.}
	\label{SS:Local3}

	Let $f: X \to S$ be a proper smooth morphism of algebraic varieties.  Then we get local systems:
	\begin{itemize}
		\item[1.] Over $\CC$: $R^i f_* \CC = \mr{Ker} (\nabla ^{an})$ where
		\[
		\nabla : H^i _{dR} (X/S) \to \Omega ^1 _{S/\CC} \otimes _{\mc{O}_S} H^i _{dR} (X/S) , 
		\quad
		H^i _{dR} (X/S) := \mbf{R}^i f_*  \Omega ^{\bullet} _{X/S} 
		\]
		is the Gauss-Manin connection on the relative deRham cohomology. Then
		$(R^i f_* \CC)_{s} = H^i (X_s ^{an}, \CC)$. This has regular singular points (Griffiths; see 
		\cite{Katz70}).
		\item[2.] $\ell$-adic: $R^i f_* \QQ _{\ell}$ is a lisse $\QQ _{\ell}$-sheaf on the \'etale topology 
		of $S$. The stalk in a geometric point $(R^i f_* \QQ _{\ell})_{\bar{s}}=
		H^i _{et}(X_s \otimes _{\kappa(s)} \kappa(\bar{s}), \QQ _{\ell})$, which has an action 
		of the Galois group $\mr{Gal}  (\kappa(\bar{s})/\kappa(s)) $.
		
	\end{itemize}

	\subsection{Example: Legendre curve.}
	\label{SS:Local4}
	Let $S= \mr{Spec}\, \ZZ [1/2, \lambda, (\lambda (1-\lambda))^{-1}]$. This is the open 
	subset  
	\[
	\mathbb{P}^1 _{\ZZ [1/2]} -\{0, 1, \infty   \} \subset \mathbb{P}^1 _{\ZZ [1/2]} .
	\]
	For each $\lambda$ we let $E_{\lambda}$ be the projective nonsingular model of the affine curve
	\[
	y ^2 = x (x-1)(x-\lambda). 
	\]
	This is an elliptic curve with origin at infinity $(x, y, z) = (0, 1, 0)$ in the projective plane. We get a proper smooth 
	morphism $f: E \to S$ defined over $\ZZ [1/2]$ with these fibers.

	The local system over $\CC$ topologically 
	is the flat bundle of $H^1 (E_{\lambda}, \CC) $ for $\lambda \in \mathbb{P}^1 (\CC) - \{0, 1, \infty \}$. This latter 
	is the sphere with 3 points removed. We have
	\[
	\pi _1 (\mathbb{P}^1 (\CC) - \{0, 1, \infty \} , \lambda _0) \cong \langle \gamma _0, \gamma _1, \gamma _{\infty} \mid
	\gamma _0 \gamma _1 \gamma _{\infty} = 1
	\rangle , 
	\]
	isomorphic to a free group on 2 generators. We can represent these by loops starting at the base-point $\lambda _0$ and
	circling once around $0, 1, \infty$, respectively. The monodromy can be represented topologically by 
	following the generators of $H^1 (E_{\lambda}, \ZZ)\sim \ZZ ^2$ as $\lambda$ moves along these loops. This can be viewed 
	by analytically continuing the period matrix of $E_{\lambda}$ relative to a basis of differential forms. The deRham 
	cohomology $H^1 _{dR}(E/S)$ is a free $\mc{O}_S$-module of rank 2.  We can take 
	as basis 
	\[
	\omega _1 = \frac{dx}{y}, \quad \omega _2 = \nabla \left ( \frac{d}{d\lambda} \right )\omega _1 = \frac{dx}{2(x - \lambda)y}. 
	\] 
	which is a basis of meromorphic 1-forms of the first and second kind modulo exact forms.

	We have the following equation in  $H^1 _{dR} (E/S)$:
	\[
	\nabla \left ( \frac{d}{d\lambda} \right )\begin{pmatrix} \omega_1 \\
		\omega _2
	\end{pmatrix} =  \frac{1}{\lambda (1-\lambda)}\begin{pmatrix}
		0 &  \lambda (1 - \lambda) \\
		1/4 &2 \lambda -1
	\end{pmatrix}  \begin{pmatrix} \omega_1 \\
		\omega _2
	\end{pmatrix} \quad \text{mod  exact forms.} 
	\]
	Note that $H^1 _{dR} (E^{an}/S^{an}) ^{\nabla = 0} = R^1 f_* \ZZ \otimes \CC $, 
	where  $R^1 f_* \ZZ$ is a  local system of free  $\ZZ$-modules of rank 2 on
	$S^{an} = \P^1 (\CC)- \{ 0, 1, \infty\}$ with a nondegenerate 
	symplectic pairing on it.  The dual local system  $R_1 f_* \ZZ $ is the homology 
	local system: $H_1 (E ^{\mr{an}}_{\lambda}, \ZZ)$. For any local horizontal section $\gamma $ of the dual  
	$R_1 f_* \ZZ $ of $R^1 f_* \ZZ $, 
	it follows easily from this that  
	the period $f(\lambda) = \int _{\gamma } dx/y$ is a solution to the differential equation

	\[
	\lambda(1-\lambda)f'' + (1-2\lambda)f' -f/4 = 0, 
	\]
	which is the differential equation satisfied by $\pFq{2}{1}{1/2&  1/2}{&1}{\lambda}$.  In fact, 
	\[
	\pFq{2}{1}{1/2& 1/2}{& 1}{\lambda}  = \frac{1}{\pi} \int _{1}^ {\infty} \frac{dx}{y}.
	\]
	We view a small loop around the interval $[1, \infty]$  as topological cycle 
	$\gamma$
	on the Riemann surface $E_{\lambda}$ which is a branched cover of
	$\mathbb{P}^1(\CC )$ of degree 2. 
	
	The monodromy group of this differential equation is projectively equivalent to the 
	principal congruence subgroup $\Gamma (2) \subset \slz$. In a suitable basis 
	(we can take so-called vanishing cycles at $1$ and $0$)
	$\gamma, \delta $ of $H^1 (E_t, \ZZ)$ the monodromy matrices are 
	\[
	T_0 = \begin{pmatrix}  1 & 2 \\ 0 & 1 \end{pmatrix}, \quad  T_1= \begin{pmatrix}  1 & 0 \\ -2 & 1 \end{pmatrix}.
	\]
	\begin{figure}[h]
		\begin{center}
			\includegraphics[width=4 cm,height=3.25cm]{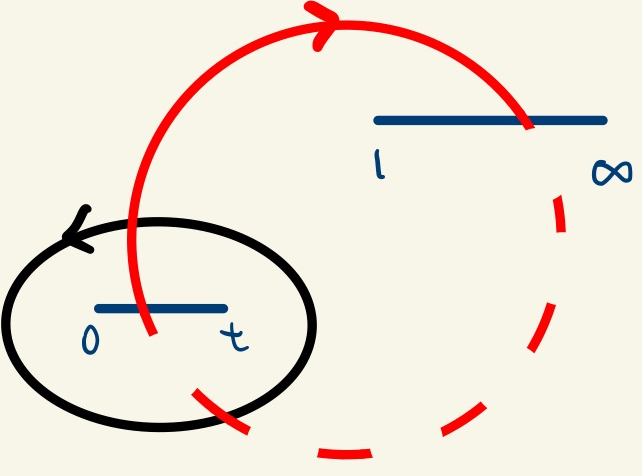} \qquad \includegraphics[width=4 cm,height=3.25cm]{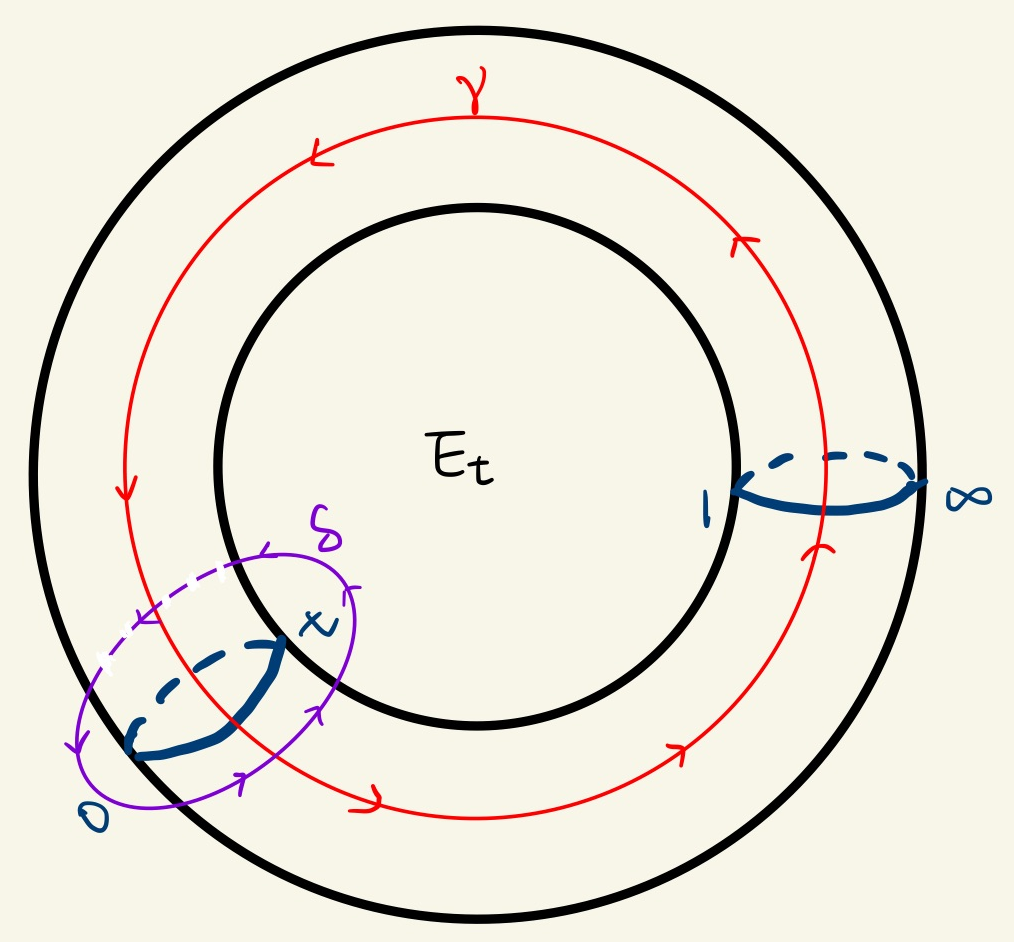}
			\caption{The Riemann surface $E_t$ is a double cover of this with branch cuts as shown.}
			\label{cycle}
		\end{center}
	\end{figure}
	
	These generate a subgroup \ft{$\G$} of index 2 inside $\Gamma (2)$ acting on the 
	upper half plane with quotient 
	\[
	\Gamma \backslash \mathfrak {H} = \Gamma (2)\backslash \mathfrak {H} 
	\overset{\sim}{\longrightarrow} \mathbb{P}^1 (\CC) - \{0, 1, \infty \}, \quad 
	\tau \mapsto \lambda(\tau).
	\]
	This map is induced by $\lambda : \mathfrak{H} \to \CC \cup \{ \infty \}$
	where $\lambda(\tau)$ is a generator of the field of modular 
	functions for $\Gamma(2)$.  The 
	inverse of this is the multivalued function 
	on $ \mathbb{P}^1 (\CC) - \{0, 1, \infty \}$
	given by the period ratio $\tau(\lambda ) = \int _{\alpha} \omega _1/\int _{\beta }\omega_1$ of two independent 
	solutions to this differential equation. The above quotient 
	space, denoted $M_2$, is the coarse moduli space of elliptic curves with a level 2 structure. 
	
	The $\ell$-adic local system $R^1 f_* \QQ _{\ell}$ gives zeta functions of the elliptic curves in the family. For instance, 
	let $\lambda \in \FF_p$ for an odd prime $p$, $\lambda \ne 0, 1$. Then this gives a geometric point 
	$\bar{\lambda} $ of the scheme $S$, and we can identify the fiber
	\[
	(R^1 f_* \QQ _{\ell} )_{\bar{\lambda}} = H^1 (E _{\lambda} \otimes \kappa(\bar{\lambda}) , \QQ _{\ell}), 
	\quad \mr{Gal}(\kappa(\bar{\lambda})/\kappa(\lambda)) = 
	\mr{Gal}(\bar{\FF} _p/\FF_p) = \langle  \mr{Frob}_p \rangle 
	\]
	which is the dual of the Tate module of the curve $E _{\lambda} $. Then
	\[
	\# E_{\lambda} (\FF_p) = p+1 - \mr{Tr} ( \mr{Frob}_p \mid (R^1 f_* \QQ _{\ell} )_{\bar{\lambda}})
	\]
	
	\begin{figure}[h]
		\begin{center}
			\includegraphics[width=5.25cm,height=5.25cm]{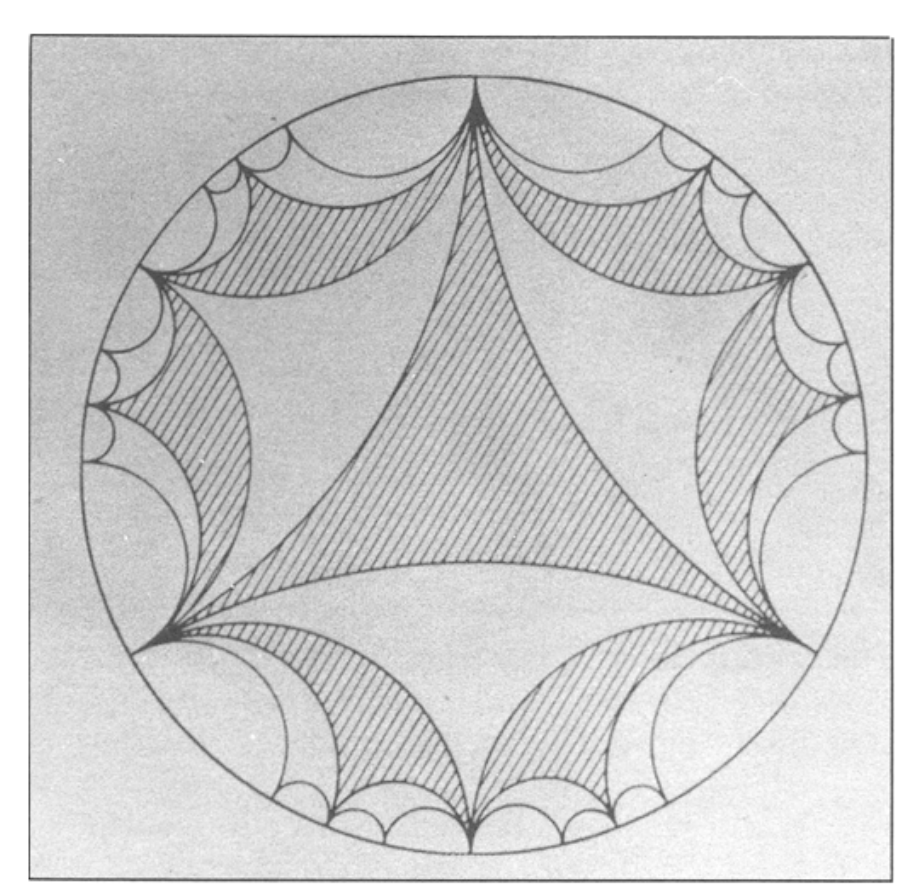}
			\caption{$(\infty, \infty, \infty)$ hyperbolic tessellation of the unit disk $\mathfrak{D}$ corresponding to $\Gamma (2)$. Image thanks to Fricke -Klein.}
		\end{center}
	\end{figure}

	\subsection{Example: Generalized Legendre curves.}
	\label{SS:Local5}
	
	For any ring $R$, we denote by $\lambda$ the coordinate on $\mathbb{A}^1 _R$ and 
	by $S$ the open set where $\lambda (1-\lambda)$ is invertible. Given any elements
	$a, b, c \in R$ we denote by $E(a, b, c)$ the free $\mc{O}_{S}$-module of rank
	2 with basis $e_0$, $e_1$,  and integrable $R$-connection
	
	\begin{align*}
		\nabla \left (\frac{d}{d\lambda}\right )(e_0) = &e_1\\
		\nabla \left (\frac{d}{d\lambda}\right )(e_1) = 
		&-\frac{(c-(a+b+1)\lambda)}{\lambda (1-\lambda)} e_1 +\frac{ab}{\lambda (1-\lambda)} e_0.
	\end{align*}
	Horizontal sections of the dual of $E(a, b, c)$ over an open set $U \subset S$ can be identified with 
	$f \in \Gamma (U, \mc{O}_U)$ which satisfy the differential equation
	\[
	{\lambda (1-\lambda)} \left (\frac{d}{d \lambda} \right )^2 f 
	+ (c -(a+b+1)\lambda)\frac{df}{d \lambda} -ab f = 0.
	\]
	These hypergeometric equations are two-dimensional factors of the the cohomology of the
	family of curves $y^N = x^a (x-1)^b(x-\lambda)^c$. In effect, the Euler integral representation
	\[
	F(\alpha, \beta; \gamma; \lambda)= \frac{\Gamma (\gamma)}{\Gamma (\beta)\Gamma(\beta - \gamma)}
	\int _1 ^{\infty} x ^{\alpha - \gamma}(x-1)^{\gamma-\beta -1}(x-\lambda)^{-\alpha}dx
	\]
	shows that the solutions to the differential equation are given by periods of those curves.  
	Given integers $N, a, b, c$ greater that zero, let $Y(N;a, b, c)_{\lambda}$ be the nonsingular projective model of the 
	affine curve in $(x, y)$-space defined by the equation  $y^N = x^a (x-1)^b(x-\lambda)^c$. 
	We consider this as a 
	family of curves 
	$$f: Y(N; a, b, c)\to  S := \mathbb{P}^1  _{R}-\{ 0, 1, \infty \}, \quad R =  R_N = \ZZ [1/N, \zeta _N].$$
	We have the Gauss-Manin connection 
	\[
	\nabla: H^1 _{dR}(Y(N;a, b, c)/S)\to \Omega ^1 _{S/R}\otimes  _{\mc{O}_S}H^1 _{dR}(Y(N;a, b, c)/S).
	\]
	The following theorem gives the structure of this, at least in the generic fiber
	$\mr{Spec}(\CC (\lambda)) \hookrightarrow U$. Choose an embedding $R \to \CC$. 
	Let
	\[
	X(N; a, b, c) = \mr{Spec} \ \CC (\lambda) [x, y, 1/y]/(y^N - x^a (x-1)^b(x-\lambda)^c).
	\]
	This is the open affine subset where $y$ is invertible. It is affine and smooth of relative 
	dimension one over $\CC(\lambda)$. The map $(x, y)\to x$ is a finite 
	\'etale covering 
	$$\pi : X(N; a, b, c) \to\mathbb{A}^1_{\CC(\lambda)}  - \{0, 1, \lambda \}. $$ 
	For any root of unity $\xi \in \mu _n$ there is an automorphism of $X(N; a, b, c) $ given by 
	$(x, y) \mapsto (x, \xi y)$. This gives the Galois group of the covering $\pi$.
	Note that 
	the $dx/y^m$ defines an element in the character eigenspace
	$ H^1 _{DR}(X(N;a, b, c)/\CC(\lambda)) ^ {\chi (-m)}$
	where $\chi (t) (\xi) = \xi ^{t}$.

	\begin{proposition}
		(\cite[6.8.6]{Katz72}) Suppose that $N$ does not divide $a, b, c, a+b+c$. Then for any
		integer $k\ge 1$ which is invertible modulo $N$ the map
		\begin{align*}
			e_0 &\mapsto \mr{the\  class\  of\ \ } \frac{dx}{y^k}\\
			e_1 &\mapsto\nabla \left (  \frac{d}{d\lambda} \right ) ( \mr{the\  class\  of\ \ }\frac{dx}{y^k})
		\end{align*}
		induces an isomorphism
		\[
		E\left( \frac{kc}{N}, \frac{k(a+b+c)}{N}-1, \frac{k(a+c)}{N}   \right ) \cong
		H^1 _{dR}(X(N;a, b, c)/\CC(\lambda)) ^ {\chi (-k)}.
		\]
	\end{proposition}
	This gives only the part of the cohomology belonging to primitive characters
	modulo $N$. For $n$ prime to $N$, we have
	$ H^1 _{dR}(X(N;a, b, c)/\CC(\lambda)) ^ {\chi (n)} = H^1 _{dR}(Y(N;a, b, c)/\CC(\lambda)) ^ {\chi (n)}$.
	The local system $R^1 f_* \CC$ on $S^{an} $ underlies a polarized variation of Hodge structures
	of weight 1. The modules with connection 
	$H^1 _{dR}(X(n;a, b, c)/\CC(\lambda)) ^ {\chi (n)}$ 
	correspond to the rank 2 local system  $(R^1 f_* \CC) ^{\chi (n)}$. The
	$(R^1 f_* \CC) ^{\chi (n)}$ does not correspond to a variation of Hodge structures unless 
	the character $\chi(n)$ is real.

	The eigenspaces  $(R^1 f_* \qlbar) ^{\chi (n)}$
	then give the $\ell$-adic realization, where $f: Y(N;a, b, c) \to S$ is as before 
	but now as schemes over $S$ for the \'etale topology. 
	
	These results are generalized to families of curves of the form
	\[
	y^N = x^a(x-1)^b (x -\lambda_1)^{c_1}....(x-\lambda_r)^{c_r}.
	\]
	in \cite{Holz86}. One gets a regular holonomic system of partial differential equations 
	in the variables $\lambda _1, ..., \lambda _r$.

	\section{Hypergeometric motives}
	\label{S:Hypmot}
	In this paper, the word {\it motive} (really: motivic sheaf) is used informally. These concepts 
	can be rigorously applied utilizing results of Arapura
	\cite{Arapura1}, \cite{Arapura2}, and Ayoub \cite{AyoubI}, \cite{AyoubII}. For our purpose, 
	the motives on a scheme $X$ will be viewed as giving 
	\begin{itemize}
		\item[1.] a constructible sheaf $\mathcal{F}_{\CC}$  of $\CC$-vector spaces on the analytic space
		$X^{an}= X(\CC)$;
		\item[2.] a constructible $\bar{\QQ}_{\ell}$- sheaf $\mathcal{F}_{\ell}$ for the \'etale topology on $X$, for a 
		fixed prime number $\ell$.
	\end{itemize}
	In this section we will describe the motivic sheaves attached to {\it hypergeometric data}.
	We utilize the formalism of hypergeometric motives as developed by Roberts and Rodriguez Villegas, see
	\cite{HypMot}. 
	A hypergeometric datum $HD$ is a pair of multi-sets  $\alpha=\{a_1,\cdots,a_n\},\beta=\{b_1=1,b_2,\cdots, b_n\}$ with $a_i,b_j\in \QQ$. 	
	A datum $HD=\left\{\alpha=\{a_1,\cdots,a_n\}, \beta=\{b_1=1, b_2, \cdots, b_n\}\right\}$ with $a_i, b_j\in\QQ^\times$ is said to be defined over $\QQ$ if the set of column vectors $\left\{\begin{pmatrix}a_1\\b_1\end{pmatrix}, \cdots, \begin{pmatrix}a_n\\b_n\end{pmatrix}\right\}$ mod $\ZZ$ is invariant  under multiplication by all $r \in (\ZZ/M\ZZ)^\times$, where $M = M(HD)$, called the level of $HD$, is the least common denominator of 
	$a_1,..., a_n, b_1,..., b_n$. 
	See Definition 1 in \cite[2.2]{LLT2}. We mainly consider primitive hypergeometric data $HD$, namely $a_i -b_j \notin \ZZ$ for all $i, j$. 
	
	Let $R_M = \ZZ[\zeta _M, 1/M]$. There is a motive $\mathcal{H}(HD)$ on 
	$\mathbf{G}_{m, R_M} = \mr{Spec}\,  R_M[x, x^{-1}]$ which 
	forms a local system of rank $n$ and is pure of weight $n-1$ on $\mathbf{G}_{m, R_M} -\{ 1 \} = \mathbb{P}^1 _{R_M} - \{ 0, 1, \infty\}$.

	\subsection{Over $\CC$.}
	\label{SS:Hypmot2}
	Let $\G(x)$ denote the Gamma function. For  $k \in \ZZ$ and $a\in \CC$, define the Pochhammer symbol $(a)_k:=\G(a+k)/\G(a)$.  
	Given $HD = \{\alpha, \beta\}$ we associate the hypergeometric function in the variable $t \in \CC$: $$\displaystyle F(\alpha,\beta;t) = \,_nF_{n-1}\left[\begin{matrix}a_1&a_2& \cdots & a_n\\ &b_2& \cdots &b_n \end{matrix} \; ; \; t \right] :=\sum_{k\ge 0} \frac{(a_1)_k\cdots(a_n)_k}{(b_1)_k\cdots(b_n)_k}t^k.$$ It satisfies the Fuchsian ordinary differential equation with three regular singular points at $0, 1, \infty$: 
	
	$$
	\left[\theta\left(\theta+b_2-1\right)\cdots \left(\theta+b_n-1\right) - t\left(\theta+a_1\right)\cdots \left(\theta+a_n\right) \right]F=0, \quad \text{where}~ \theta:=t\frac{d}{dt}.
	$$
	The local solutions of this equation form a rank-$n$ local system, denoted by $\mathcal{H}(HD)_{\CC}$. 
	The local exponents of the above hypergeometric differential equation   are 
	\begin{equation}
		\begin{split}\label{eq:indicial}
			0,1-b_2,\cdots,1-b_n  & \quad \text{ at } t=0,\\
			a_1,\,a_2,\,\cdots,\,a_n &  \quad \text{ at } t=\infty,\\
			0,1,2,\cdots,n-2,\gamma& \quad \text{ at } t=1,
		\end{split}
	\end{equation} respectively, where \begin{equation}\label{eq:gamma}
		\gamma=-1+\sum_{j=1}^n b_j - \sum_{j=1}^n a_j. 			\end{equation}  
	See \cite{Beukers-Heckman, Slater}  for more details.

	\subsection{Over a finite field.}
	\label{SSS:Hypmot3}

	Let $\FF_q$ be a finite field of odd characteristic and use $ \widehat{\FF_q^\times}$ to denote the group of multiplicative characters of $\FF_q^\times$.   
	We use $\eps_q$ or simply $\eps$ to denote the trivial character, and  $\phi_q$ or  $\phi$ to denote the quadratic character. For any $A\in \widehat{\FF_q^\times}$, use $\overline A$ to denote its inverse, and extend $A$ to $\FF_q$ by setting 
	$A(0) = 0$.  For any characters  $A_i$, $B_i \in \widehat{\FF_q^\times}$, $i=1,\cdots,n$ with $B_1=\eps$ in $\widehat{\FF_q^\times}$, the $_{n}\P_{n-1}$-function (cf. \cite{Win3X}) is defined as follows: 
	
	\begin{multline}\label{eq:PP}
		\pPPq{n}{n-1}{A_1& A_2&\cdots &A_{n}}{&  B_2&\cdots &B_{n}}{\lambda;q}\\
		: = \prod_{i=2}^{n} \left(-A_iB_i(-1)\right)\cdot  \left (\frac{1}{q-1}
		\sum_{\chi\in \widehat{\FF_q^\times}}\CCe{A_1\chi}{\chi} \CCe{A_2\chi}{B_2\chi}\cdots \CCe{A_{n}\chi}{B_{n}\chi}\chi(\l)
		+\delta(\lambda) \prod_{i=2}^{n} \CCe{A_{i}}{B_i}\right),
	\end{multline}

	where
	$$   \displaystyle\CCe AB :=-B(-1)J(A,\ol B), \quad J(A, B):=\sum_{t\in \FF_q}A(t)B(1-t), \quad   \mr{ and } \quad \delta(\l):=
	\begin{cases}
		1,& \mbox{ if } \l=0, \\
		0,& \mbox{ otherwise. } \end{cases}
	$$
	It will be useful for us to have some notation for characters of finite fields. First we recall the 
	$M$th power residue symbol. Let $R_M = \ZZ[1/M, \zeta _M]$. For every prime ideal 
	$\mfr{p} \subset R_M$, $R_M/\mfr{p}$ is a finite field with $q = q(\mfr{p}) = N\mfr{p}$ elements. 
	Recall that $\FF _q ^{\times}$ is a cyclic group with $q-1$ elements, and since $\zeta _M$ maps 
	to an element of order $M$ modulo $\mfr{p}$, we have $q \equiv 1$ mod $M$. For any 
	$x \in R_M$ prime to $\mfr{p}$, there is a unique $M$th root of unity with the property
	\[
	x^{(N\mfr{p}-1)/M} \equiv 
	\left (\frac{x}{\mfr{p}}\right)_M x \text{\ mod\ } \mfr{p}, \quad \left (\frac{x}{\mfr{p}}\right) _M\in  \mu _M (\bar{\QQ}).
	\]
	The map $x \mapsto \left (\frac{x}{\mfr{p}}\right) _M$ induces a character 
	$(R_M /\mfr{p})^{\times } \to \mu _M (\bar{\QQ})$. Moreover, every such character 
	is of the form $\left (\frac{\cdot}{\mfr{p}}\right) ^i _M$ for some $i \in \ZZ/M$. Given 
	a rational number $a$ with denominator $M$, we can identify it with 
	an element of $(1/M)\ZZ/\ZZ \ = \ZZ/M$. Then we define a character 
	$\iota _{\mfr{p}}(a) (\cdot)= \left (\frac{\cdot}{\mfr{p}}\right ) ^a_M: (R_M /\mfr{p})^{\times } \to \mu _M (\bar{\QQ})$. 
	
	We choose an isomorphism 
	of all roots of unity $\QQ (\zeta _{\infty})\subset \bar{\QQ}_{\ell}$ with 
	the roots of unity $\QQ (\zeta _{\infty})\subset \CC$.

	\begin{theorem}[Katz \cite{Katz, Katz09}]\label{thm:Katz}Let  $\ell$ be a prime. Given a primitive hypergeometric datum HD = $\{\alpha, \beta\}$ consisting  of $\alpha=\{a_1,\cdots,a_n\}$, $\beta=\{1,b_2,\cdots,b_n\}$ with $a_i, b_j \in \QQ^\times$ and $M:=M(HD) = \mr{lcd}(\alpha \cup \beta)$. There exists a constructible $\bar{\QQ}_{\ell}$-sheaf for the \'etale topology on 
		$\mathbf{G}_{m, R_{M}}[1/\ell] $, denoted
		$\mathcal{H}^P(HD)_{\ell}$, with the following properties:
		
		\begin{itemize}
			\item[1.] Let $\lambda \in  | \mathbf{G}_{m, R_{M}}[1/\ell]  |$ be a closed point and 
			$\bar{\lambda}$ a geometric point localized in $\lambda$. The residue field 
			$k(\lambda)$ is a finite extension of the field $R_{M}/\mfr{p} := k({\mfr{p}})  = \FF_{q(\mfr{p})}  = \FF_q$. 
			We let $\mr{deg}(\lambda)$ be the degree of this extension. Then
			\[
			\Tr\left (  \Frob _{\lambda}  \mid       \mathcal{H}^P(HD)_{\ell, \bar{\lambda}}     \right )
			=  \pPPq{n}{n-1}{\iota _{\lambda}(a_1)&\iota _{\lambda}(a_2) &\cdots & \iota _{\lambda}(a_n)  }   
			{&  \iota _{\lambda}(b_2) &\cdots &\iota _{\lambda}(b_n)    }{1/\lambda;q(\lambda)}
			\]
			where $q(\lambda) = q^{\mr{deg}(\lambda)}$ and $\iota_{\lambda}(c)$ is the character
			\[
			x \mapsto  \iota _{\mfr{p}} (c)(  N _{k(\lambda)/k(\mfr{p})}    (x)), \quad
			k(\lambda) ^{\times}\to \mu _M(\bar{\QQ}) \subset \CC ^{\times}.
			\]
			
			\item[2.] When $\l\neq 1$,  the stalk $(\mathcal {H}^P(HD)_\ell)_{\bar{\l}}$ has dimension
			$n$ and all roots of the characteristic polynomial of $\Frob_{\lambda}$  are algebraic numbers and have the same absolute value $q(\lambda)^{(n-1)/2}$ under all archimedean embeddings. 
			
		\end{itemize}

	\end{theorem}

	There is a variant of these character sums due to Beukers/Cohen/Mellit
	which has the important property that in many cases
	we can define these sheaves slightly more generally. For instance, 
	suppose that the data $HD$ is defined over $\QQ$, then there is a sheaf
	$\mathcal{H} ^{BCM}(HD)_{\ell}$ on $\mr{Spec}(\ZZ[1/M\ell, \lambda, \lambda^{-1}]  $
	such that 
	\[
	\mathcal {H} ^{BCM}(HD)_{\ell}\otimes  _{\ZZ[1/M\ell]}\ZZ[1/M\ell, \zeta _M] = \mathcal {H} ^P(HD)_{\ell}.
	\]
	The Frobenius traces of $\mathcal {H} ^{BCM}(HD)_{\ell}$ are $\QQ$-valued.

	\subsection{Transformations}
	\label{SS:Hypmot4}
	There is a huge number of identities relating different hypergeometric functions. These have 
	parallel versions over $\CC$ and over finite fields. For instance:
	
	The Clausen formula (\cite{AAR} or \cite[Eqn. (39)]{LLT2}):   
	\begin{equation}\label{eq:clausen1}
		{(1-t)^{-\f12}} \pFq32{\f12&a-b+\f12&b-a+\f12}{&a+b+\f12&\frac32-a-b}{t}\\=  \pFq21{a&b}{&a+b+\f12}{t}\pFq21{1-a&1-b}{&\frac32-a-b}{t}
	\end{equation}when both sides are convergent. The finite field analog

	\begin{theorem}[Evans-Greene \cite{Evans-Greene}]\label{thm:E-G} 
		For a fixed finite field $\FF_q$, let $\phi$ be the quadratic character.
		Assume $\eta,K\in \widehat{\FF_q^\times}$ such that none of $\eta, K\phi,\eta K, \eta \ol K$ is trivial.  
		When $t=1$, we have
		\begin{multline*}
			\pPPq32{\phi&\eta&\ol \eta}{&K&\ol K}{1} 
			=\begin{cases}
				0,& \mbox{ if $\eta K$ is not a square in $\widehat{\FF_q^\times}$,}\\
				\frac{J(\eta K,\ol \eta K)}{J(\phi, \ol K)}\(J(S \ol K, \phi\ol S)^2 + J(\phi S\ol K, \ol S)^2\),& \mbox{ if $\eta K=S^2$  in $\widehat{\FF_q^\times}$}.
			\end{cases}
		\end{multline*}
		
		When $t\neq 0,1$, suppose $\eta K = S^2$ for some  $S\in \widehat{\FF_q^\times}$, we have 
		\begin{align*}
			\phi(1-t)  \pPPq32{\phi&\eta&\ol \eta}{&K&\ol K}{t}
			&=  \pPPq21{\phi K \ol S& S}{&K}t  \pPPq21{\phi \ol K S & \ol S}{& \ol K}t -  q.
		\end{align*} 
	\end{theorem} 
	
	For a finite field analog of the Whipple $\phantom{}_7 F_6$-formula see \cite{LLT2}.
	For an application of a cubic transformation formula for an Appell-Lauricella hypergeometric
	function over a finite field, see \cite{FST}.

	\section{Automorphic motives}
	\label{S:autmot}
	There are motives $V^n (\Gamma)$ on modular curves $X_{\Gamma}$ such that the 
	cohomology $H^1 (X_{\Gamma}, V^n (\Gamma) )$ is related to cusp forms on $\Gamma$ 
	of weight $n+2$. This is the content of Eichler-Shimura theory. There are two main aspects: the geometric description of modular forms, 
	and the congruence formula
	relating the trace of Hecke operators to the trace of Frobenius.

	\subsection{Over $\CC$.}
	\label{SSS:autmot1}
	The main result here due to Eichler \cite{Eichler57} and Shimura \cite{Shimura59} is an isomorphism 
	\[
	S_{n+2}(\Gamma) \overset{\sim}{\longrightarrow} H^1 _{par} (\Gamma, V^n _{\RR})
	\]
	between the space of cusp forms of weight $n+2$ for $\Gamma$ and the parabolic cohomology group, 
	where $V^n _{\RR}$ is the $n$-th symmetric power of the standard representation of 
	$\slr$ on $\RR ^2$. Parabolic cohomology is treated in detail in Shimura's book, \cite{Shimura-book}.
	
	Here is a brief outline: let $\begin{pmatrix}u \\  v\end{pmatrix} \in \CC^2$. For any $n \ge 0$, let
	$\begin{pmatrix}u \\  v\end{pmatrix} ^n\in \CC^{n+1}$ be the vector whose components 
	are $u^n, u^{n-1}v, ..., v^n$ ($= 1$ if $n=0$.) We define a representation $\rho_n$ of $\mr{GL}(2, \CC)$ of dimension 
	$n+1$ by the rule:
	\[
	\rho_n(\alpha)\begin{pmatrix}u \\  v\end{pmatrix} ^n 
	:= \left (   \alpha \begin{pmatrix}  u \\ v  \end{pmatrix}  \right )^n,  \quad
	\alpha \in \mr{GL}(2, \CC). 
	\]
	That is, $\rho_n (\alpha) = \mr{Sym}^n (\alpha)$. One checks, for $\alpha \in \mr{GL}(2, \RR) ^+$:
	\[
	\begin{pmatrix}\alpha(z) \\  1\end{pmatrix} ^n = j(\alpha, z)^{-n} \rho_n(\alpha) \begin{pmatrix}z \\  1\end{pmatrix} ^n , \quad
	\text{where \ }
	z \in \mathfrak{H}, \alpha (z) = \frac{a z+b}{cz+d}, \alpha = \begin{pmatrix} a & b \\c& d \end{pmatrix}, 
	j(\alpha, z ) = (c z +d).
	\]
	We have
	\[
	j(\alpha \beta, z) = j(\alpha,  \beta( z))j( \beta, z).
	\]
	Let $\Gamma \subset \mr{SL}(2, \RR)$ be a Fuchsian subgroup of the first kind and $n\ge 0$
	an integer. If $-1 \in \Gamma$ we take $n$ to be even. If $f \in S_{n+2}$ is a cusp-form 
	of weight $n+2$ we define a vector-valued differential form
	\[
	\mathfrak{d}(f)= \begin{pmatrix} f(z)z^n dz \\f(z)  z^{n-1} dz \\ \vdots \\f(z)dz\end{pmatrix}.
	\]
	We get
	\[
	\mfr{d}(f )\circ \alpha  = \rho_n(\alpha )\mfr{d}(f),\quad
	\mr{Re (}\mfr{d}(f ))\circ \alpha  = \rho_n(\alpha ) \mr{Re} ( \mfr{d}(f)),  \quad \alpha\in\Gamma,
	\]
	the second equality holding because the representation $\rho_n$ is real. Now fix any point 
	$z_0 \in \mfr{H}^*$ and define for $f \in S_{n+2}(\Gamma)$ the indefinite integral
	\[
	F(z) = \int _{z_0}^z \mfr{d}(f) + v, \quad z\in \mfr{H}^*, v \in \CC^{n+1}, 
	\]
	for any path in $\mfr{H}$ connecting $z_0$ and $z$. If they are both 
	cusps, we can take a geodesic arc in $\mfr{H}$ connecting them.
	Then
	\[
	F(\alpha (z)) = \rho_n(\alpha )F(z) + t(\alpha)
	\]
	where $\alpha \mapsto t(\alpha)$ is a parabolic 1-cocyle of $\Gamma$ with values
	in the representation $V^n _{\CC}$. This means
	\begin{itemize}
		\item[1.] $t(\alpha \beta) = t(\alpha)+ \rho_n(\alpha)t(\beta)$.\
		\item[2.] For each parabolic element $\sigma \in \Gamma$ there is a vector $w \in V^n_{\CC}$
		such that $t(\sigma ) = w -\rho_n(\sigma)w$. 
		The vector $w$ may depend on the parabolic element $\sigma$.
	\end{itemize}
	We define $H^1 _{par}(\Gamma, V^n_{\CC})$ as the quotient of all parabolic 1-cocyles (property 1 and 2)
	modulo all coboundaries, that is,  
	there is a $w$ such that $t(\alpha ) = w -\rho_n(\alpha)w $ for all $\alpha \in \Gamma$. Similarly, taking real parts defining
	\[
	F_{\RR}(z) = \int _{z_0}^z \mr{Re}(\mfr{d}(f)) + v, \quad z\in \mfr{H}^*, v \in \RR^{n+1},  
	\]
	the corresponding cocycle $u(\alpha)$ takes values in $H^1 _{par} (\Gamma, V^n _{\RR})$ and depends
	only on $f$ and not on the auxiliary choices of $z_0$ and $v$. 
	Shimura proved that $f \mapsto u(\alpha)$ defines an isomorphism of real vector spaces. 
	$S_{n +2}(\Gamma)\cong H^1 _{par} (\Gamma, V^n _{\RR})$. This implies that 
	$\dim H^1 _{par} (\Gamma, V^n _{\CC}) = 2 \dim S_{n+2}(\Gamma)$. When $\Gamma$ is an 
	arithmetically defined group, there are Hecke operators defined and the Eichler-Shimura isomorphism
	respects this action. Moreover in that case there is a lattice $V_{\ZZ}^n \subset V_{\RR}^n $ stabilized 
	by $\Gamma$ and we have parabolic cohomology groups $H^1 _{par}(\Gamma, V_{\ZZ} ^n)$. One can 
	define a complex structure and a polarization such that the torus $H^1 _{par}(\Gamma, V_{\RR} ^n)/H^1 _{par}(\Gamma, V_{\ZZ} ^n)$
	has the structure of an abelian variety $A(\Gamma, n)$. 
	Finally, if $f\in S_{n+2}(\Gamma)$ is a Hecke eigenform the periods of $f$ are related to the values of the 
	$L$-function $L(f, s)$ at integer points in the critical strip for this Dirichlet series. See Manin's papers, 
	\cite{Manin72}, \cite{Manin73} 
	for the connection to modular symbols, $p$-adic Hecke series and more.  See   Greenberg and Voight \cite{Voight2011} for 
	algorithms for computing Hecke eigenvalues in parabolic cohomology. See Stiller's monograph \cite{Stiller84} for more
	information on the connection to differential equations and special values of $L$-functions. See the papers of Zagier \cite{Zagier-top-diff} and Kontsevich/Zagier
	\cite{K-Z} for more information on differential equations and periods.

	Here are examples from \cite{Zagier-top-diff, LLT2}:
	$$
	\pFq43{\f12&\f12&\f12&\f12}{&1&1&1}{1}=\frac{16}{\pi^2}L(\eta(2\tau)^4\eta(4\tau)^4,2),
	$$
	$$
	\pFq65{\f12&\f12&\f12&\f12&\f12&\f12}{&1&1&1&1&1}{1} =16\oint_{|t_2(\tau)|=1} \tau^2 f_{8.6.a.a}\left (\frac{\tau}2\right)d\tau, 
	$$
	where 
	$ t_2(\tau)=-64\frac{\eta(2\tau)^{24}}{\eta(\tau)^{24}}$ is a Hauptmodul for $\Gamma_0(2)$, and $f_{8.6.a.a}$ is a normalized cuspidal newform expressed in LMFDB label.

	Deligne \cite{Deligne71} reinterpreted this in the following manner. Assume that $\Gamma$ has no torsion. 
	First, we have 
	an isomorphism
	\[
	S_{k+2}(\Gamma) \cong H^0 (X _{\Gamma},  \omega^k \otimes \Omega^1 _{X _{\Gamma},} )
	\]
	for certain coherent sheaves $\omega ^k$ on the compact Riemann surface $X ^{an}_{\Gamma}$. These 
	sheaves have the property that 
	\[
	\omega ^2 \cong  \Omega^1 _{X _{\Gamma}} (\mr {log} S_{\Gamma}), \quad S_{\Gamma} =
	\text{the set of cusps of\ } X _{\Gamma}.
	\]
	In fact, over the complement of the cusps, the sheaf $\omega ^k$ is the holomorphic
	line bundle attached to the cocycle $j(\alpha, z)^k$. 
	There are constructible sheaves of $\CC$-vector spaces $V^k (\Gamma)_{\CC}$
	on $X_{\Gamma}  $ for integers $k\ge 1$ with the following properties:
	\begin{itemize}
		\item[1.] Over the open subset $X_{\Gamma} ^{\circ}$  
		these  are local systems of rank 
		$k+1$. They are thus solution sheaves to a system of regular singular differential equations
		on  $X_{\Gamma} ^{\circ}$. 
		For $\Gamma$ considered in this paper, these will be hypergeometric  local systems when $k=1$ or 2.
		
		\item[2.] The sheaves $V^k (\Gamma)_{\CC}$ on $X_{\Gamma} $ are extensions of the sheaves on $X_{\Gamma} ^{\circ}$ defined above: 
		\[
		V^k (\Gamma) _{\CC}  = \iota_*   (V^k (\Gamma) _{\CC} \mid X_{\Gamma} ^{\circ}) \quad \text{for~the~inclusion~map}~
		\iota: X_{\Gamma} ^{\circ}  \to X_{\Gamma}.
		\] 
		\item[3.] We have:
		\begin{equation*}
			\label{E:mot4}
			H^1 _{par} (\Gamma , \mr{Sym}^k (\CC ^2))  \cong 
			H^1 (X_{\Gamma} , V^k (\Gamma)_{\CC} )  \cong
			S_{k+2} (\Gamma) \oplus \overline{ S_{k+2} (\Gamma)}, 
		\end{equation*}
		which is a Hodge decomposition of type $(k+1, 0), (0, k+1)$.  
		The left-hand side is parabolic cohomology. 
		Moreover
		\[
		H^1 (X_{\Gamma} , V^k (\Gamma)_{\CC} ) \cong \mr{Im}(
		H^1 _c (X_{\Gamma}^{\circ} , V^k (\Gamma)_{\CC} )\to 
		H^1  (X_{\Gamma}^{\circ} , V^k (\Gamma)_{\CC} 
		))
		\]
		the image of the compactly supported cohomology.
		The Hecke operators act on the spaces as geometric correspondences and these
		isomorphisms are equivariant for the Hecke actions. 
	\end{itemize}
	When $\G\subset  \mr{SL}_2(\RR)$ is an arbitrary Fuchsian subgroup of the first kind, proofs of these theorems can be found in 
	\cite{BayerNeu81}. When $\Gamma$ is torsion-free, these are defined as follows. 
	The group $\Gamma$ acts on 
	$\mathfrak{H}^*\times \mr{Sym}^k (\CC^2)$, on the first factor by linear fractional transformations, 
	on the second factor via $\varrho _k := \mr{Sym}^k (\varrho_1)$,  
	where $\varrho_1$
	arises from the inclusion $\G\subset  \mr{SL}_2(\RR)$ and the canonical action of 
	$  \mr{SL}_2(\RR) $ on $\CC^2$. Then 
	\begin{equation}\label{eq:V-CC}
		V^k (\Gamma)_{\CC}:= p_* ^{\G} (\underline{\mr{Sym}^k (\CC ^2)}) = \G\backslash 
		\mathfrak{H}^*\times \mr{Sym}^k (\CC ^2).
	\end{equation}
	When $\Gamma$ has torsion, one takes a normal subgroup $\Gamma' \subset \Gamma$ 
	of finite index and 
	defines 
	\[
	V^k (\Gamma) := p_* ^{\Gamma/\Gamma'}(V^k (\Gamma')). 
	\]
	This works with certain restrictions, namely if $-1 \in \Gamma$ this becomes identically 0
	when $k$ is odd, so one must restrict to even $k$ (this is no limitation anyway since
	there are no cusp forms of odd weight if  $-1 \in \Gamma$).

	\subsection{$\ell$-adic.}
	\label{SS:autmot2}
	Original papers: \cite{Eichler54}, \cite{Shimura58}, \cite{Shimura61}. It is preferable 
	to replace the group $\mathrm{SL}_2$ by  $\mathrm{GL}_2$ and work with the adeles. 
	Let $G$ be the algebraic group $\textsf{GL}_2$ or $\textsf B^\times$ according as $\G$ non-cocompact or cocompact arising from an indefinite quaternion division algebra $B$ defined over $\QQ$. 
	In both cases, $G(\RR) =$ GL$_2(\RR)$  which has two connected components $G(\RR)_{\pm}$ consisting of elements 
	with positive and negative determinants, respectively.  
	Moreover, $G(\RR)_+$ acts transitively on the upper half plane 
	$\mathfrak H$ via fractional linear transformations so that we may identify $\mathfrak H$  with $G(\RR)_+/\rm{SO}_2 Z(\RR)_{+}$, where $Z$ denotes the center of $G$.  Thus  $X:= G(\RR)/\rm{SO}_2Z(\RR)$ may be identified with $\mathfrak H^{\pm}$, the disjoint union of upper and lower half plane. The pair $(G, X)$ satisfies the axioms of Deligne 
	for Shimura varieties, \cite{Deligne79}.  Let $G(\QQ)_+ =G(\QQ)\cap G(\RR)_+$. Write $\mathbf A$ for the ring of adeles over $\QQ$, and $\mathbf A^f$ the subring of finite adeles.  Let $U$ be a compact-open subgroup of $G(\mathbf A_f)$.  
	Define
	\[
	Y_U  ^{ an}:= \phantom{}_U M_{\CC}(G, X) = G(\QQ)_{+}\backslash (X^{+}\times G( \mathbf A^f  ))/U . 
	\]
	This is a finite disjoint union, indexed by $ S_{U} = G(\QQ)_{+}\backslash G( \mathbf A^f  )/U$, 
	of quotients $\Gamma _g\backslash X^+$, for subgroups 
	$\Gamma _g \subset G^{\mr{ad}}(\RR)^+$ which are the images of 
	$\Gamma ' _g = g U g^{-1}\cap G(\QQ)_+$ under $G(\QQ)_+ \to G(\RR)_+\to  G^{\mr{ad}}(\RR)^+$, 
	where $G^{\mr{ad}} = G/Z$.  Shimura's theory of canonical models gives a curve $Y_U$ defined over $\QQ$ whose complex points are these. This curve is irreducible, but not absolutely irreducible. The irreducible components 
	of $Y_U \otimes _{\QQ}\bar{\QQ}$ are indexed 
	by the same finite set as above, and the corresponding analytic spaces are the  $\Gamma _g\backslash X^+$. These
	components are defined over specific abelian number fields $k_{\Gamma _g}$. 
	
	We also have the spaces of cusp forms $S_k (U)$. We can describe this as 
	\[
	\bigoplus _{g\in S_U} S_{k} (\Gamma _g\backslash X^+), 
	\]
	but more intrinsically as the set of functions $\phi : G(\mathbf A)\to \CC$ such that 
	\begin{itemize}
		\item[1.] For all $\gamma \in G(\QQ)$, and $u \in U$
		\[
		\phi (\gamma g u ) = \phi(g). 
		\]
		\item [2.] For any $g_f \in G(\mathbf A ^f)$, the function $g_{\infty}\mapsto \phi(g_{\infty}g_f)$ 
		is invariant under $Z(\RR)  \cap G(\RR)_+$, smooth
		on $G(\RR)_+$ and satisfies 
		\[
		\phi(g k_{\theta}) = \phi(g) (e^{i \theta })^k, \quad k_{\theta} = \begin{pmatrix} \cos (\theta) & \sin (\theta)\\
			-\sin (\theta) & \cos(\theta)
		\end{pmatrix}
		\]
		and 
		\[
		X_{-} \phi = 0, 
		\]
		where 
		\[
		X_{\pm} = \frac{1}{2} \begin{pmatrix}
			1 & \pm i \\ \pm i & -1
		\end{pmatrix} \in \mathfrak {g}_{\CC} = \mr{Lie} (G(\CC)).
		\]
		\item [3.] If $G = {\sf GL}_2$, then for all $g \in G(\mathbf A)$, 
		\[
		\int _{\QQ \backslash \mathbf A} \phi \left(  \begin{pmatrix}  1 & x\\0 & 1 \end{pmatrix} g   \right ) dx = 0.
		\]
	\end{itemize}
	The relation between this adelic viewpoint and classical cusp forms is explained in several places, see \cite{JL}.
	The above $\phi$ corresponds to functions $f(\tau )$ which are modular with respect to 
	subgroups $\Gamma _g$ as above.  The second condition expresses the holomorphy of $f(\tau)$, 
	and that it has weight $k$. The third condition expresses the vanishing 
	of the zeroth Fourier coefficients at the cusps 
	and is only meaningful for us in the case $G = \textsf{GL}_2$.  
	
	The main theorem is 
	\begin{theorem}
		\label{T:deligneohta}
		Let $U$ be an open compact subgroup as above, stable under the canonical involution. 
		Let $d = \dim _{\CC}S_k(U)$. For each rational prime $\ell$ there exists an $\ell$-adic representation
		\[
		\psi _{U, k}: \mr{Gal} (\bar{\QQ}/\QQ) \longrightarrow
		\mr{GL}(2d, \QQ _{\ell})
		\] 
		which is unramified outside a finite set $S$ of prime numbers. For $p \notin S$	we have
		\[
		\det \left (1 -   \psi _{U, k} (\mr{Frob}_p) x \right ) 
		= \det\left  ( 1 - T(p) x + p T(p, p)x^2\mid S_k(U)
		\right )
		\]
		where $T(p)$ and $T(p, p)$ are the standard Hecke operators. 
	\end{theorem}	
	For $G = {\sf GL}_2$ this is due to Deligne; for $G = {\sf B}$ this is main result of
	Ohta, stated only for quaternion algebras over $\QQ$. For even weights $2k$ and for 
	$U$ which are sufficiently small, the $\ell$-adic representation is on the 
	$\QQ _{\ell}$-vector space
	\[
	H^1 _{et}(X_U \otimes _{\QQ } \bar{\QQ},    V^{2k} (U)_{\ell})
	\]	
	where $V^{2k} (U)_{\ell}$ are $\ell$-adic analogs of the local systems denoted
	$V^{2k} (\Gamma)_{\CC}$ in the previous section. They have similar properties: 
	they are constructible $\QQ _{\ell}$- sheaves, lisse of rank $2k+1$ on the complement of cusps
	and elliptic points. The situation for odd weights $k$ is more complicated in general.  
	In the quaternion case, there are sheaves $V^k (U)_{\ell}$ but their 
	rank is twice the corresponding $V^{k} (\Gamma)_{\CC}$. Moreover, one 
	must exclude cases where $-1 \in U$. 	
	
	The curves and sheaves all have good reduction modulo $p$ for $p \notin S$. Using the same 
	letters to denote them we get as a corollary that 
	\begin{equation}
		\label{E:ES}
		\mr{Tr}\,   (\mr{Frob}_p \mid  H^1 _{et}(X_U \otimes _{\FF _p } \bar{\FF} _p,    V^{2k} (U)_{\ell})) = \mr{Tr}\, 
		(T(p) \mid S_{2k}(U)).
	\end{equation}
	Since $ H^i _{et}(X_U \otimes _{\FF _p } \bar{\FF} _p,    V^{2k} (U)_{\ell})=0$ for all $i \ne 1$, we can use the Grothendieck-Lefschetz trace formula to calculate the left-hand side. 
	This will be explained in section \ref{S:main}.

	\section{Shimura curves.}
	\label{S:D=6}
	\subsection{Shimura curves in general}
	\label{SS:D=6.1}
	Let $F$ be a totally real number field of degree $g$ over $\QQ$. Let $B$ be an indefinite quaternion 
	algebra over $F$ such that $B \otimes _{\QQ}\RR \cong M_2 (\RR) \times \mathbb{H}^{g-1}$
	where  $\mathbb{H}$ is Hamilton's quaternions. We let $O_B$ be a maximal order (unique up 
	to conjugation) and $O_B ^1$ the set of elements of reduced norm 1. 
	If $\Gamma$ is commensurable with $O_B ^1$  then we have a compact Riemann surface arising from $\Gamma$, say $X_{\Gamma}$.   Shimura
	proved that this has a canonical model over a number field $k_{\Gamma}$. Unless $g=1$ however, this has no simple
	moduli interpretation. Nonetheless this belongs to the general theory of Shimura varieties (but not of PEL type).   
	When $F = \QQ$ Kuga and Shimura studied these and in particular the zeta functions of the fiber spaces of abelian varieties 
	over these curves. See \cite{Kuga3}, \cite{KS2}, \cite{Shimura68}. This was extended by Ohta in a series of papers, 
	\cite{Ohta82}, \cite{Ohta81}, \cite{Ohta83}. Langlands also established these results by his methods. 
	See \cite{Langantwerp},  \cite{Langlands79}. 
	These methods are now part of the standard toolbox
	in the arithmetic of Shimura varieties. In our paper we follow Kuga, Shimura and Ohta.	
	
	In Takeuchi's list of arithmetic triangle groups is the famous $(2, 3, 7)$-tesselation. 
	As Fricke discovered, this is related to a quaternion algebra over the cubic field 
	$F = \QQ (\zeta _7 + \bar{\zeta}_7)$. He also noted that the corresponding curve 
	is the Klein quartic
	\[
	x^3 y + y^3z + z^3x= 0. 
	\]
	For a modern exposition, see Elkies' article, \cite{Elkies99}.
	\subsection{Shimura curves over $\QQ$.}
	\label{SS:D=6.2}
	References for this section: \cite{Elkies}, \cite{Jordan}, \cite{Rotger04}, \cite{Voight-book}.
	Let $B$ be an indefinite quaternion division algebra over $\QQ$ of discriminant $D$. It is known 
	that $B$ ramifies at an even number primes $\Sigma$. We let $O_B = O$ be a maximal order 
	in $B$ (unique up to conjugation). Define
	\[
	\Gamma  = \Gamma _B := O^1/{\pm 1}, \quad 
	\Gamma ^* =  \Gamma^* _B := \left \{ x \in B^*/\QQ^*\mid x  O =  Ox, \mr{nr}(x) > 0
	\right \}. 
	\]
	It is known that $\Gamma $ is a normal subgroup of  $\Gamma ^*  $ with quotient an elementary
	abelian 2-group with $\#\Sigma$ generators. It is also known that the elements of $\Gamma ^*  $
	are the classes mod $\QQ^*$ of elements of ${O}$ with reduced norm 
	$\prod_{p \in T}p$ for some possibly empty subset $T \subset \Sigma$. For any subset $T \subset \Sigma$
	we let $\Gamma \subset \Gamma _T \subset \Gamma ^*$ be the subgroup corresponding to 
	$T$. 
	
	Fix an isomorphism $\theta : B \otimes _{\QQ}\RR \cong M_2 (\RR)$. Then the groups 
	$\Gamma , \Gamma ^* $ become isomorphic to subgroups of 
	$\mr{PSL}_2( \RR) = \mr {Aut}(\mathfrak{H})$. For any subset $T \subset \Sigma$
	we have quotient Riemann surfaces
	$X_{\Gamma , T} = X_{B , T}:= \Gamma _T \backslash \mathfrak{H}$. These are compact (no cusps), but there
	are elliptic points. Shimura's theory of canonical models shows that these are the sets of complex
	points of an algebraic curve, denoted by the same symbol if no confusion is possible, 
	defined over a specific number field $k_{B, T}$. These have interpretations as coarse moduli 
	spaces of polarized abelian surfaces, as we briefly recall. 
	
	To define a polarization we choose an element $\mu \in O$ with 
	$\mu ^2 = -D$. This defines an anti-involution $a \mapsto a'$ of $B$  
	by the rule $a ' = \mu ^{-1}a^* \mu$ where $a^*$ is the canonical involution:
	$\mr{nr}(a) = a a^*, \mr{tr}(a) = a + a^*$. This is a positive involution in that the quadratic
	form $a \mapsto \mr{tr}(a'a)$ is positive definite. We describe the moduli space 
	in the case $T = \emptyset$, i.e., for $\Gamma$. Namely $X_{B} = X_{\Gamma}$ is the coarse 
	moduli space of triples $(A, \rho, \iota)$ where $A$ is an abelian surface, 
	$\rho$ is a principal polarization and $\iota: O\to \mr{End} (A)$ is an embedding
	such that the Rosati involution defined by $\rho$ on $\iota (O)$ is the involution 
	defined by $\mu$. Concretely, for each $z \in \mathfrak{H}$ we can define a 
	triple $(A_z, \rho _z, \iota _z)$ as follows. 
	\[
	A_z = \CC^2 / \Lambda _z , \quad      \Lambda _z :=   \theta (O) v_z, \quad 
	v_z = \begin{pmatrix}   z \\ 1   \end{pmatrix}, \quad  \theta (O)\subset M_2 (\RR ).
	\]
	Then $E_z : \Lambda _z \times \Lambda _z \to \ZZ$ given by 
	\[
	E_z (\theta (\lambda _1)v_z,\theta (\lambda _2)v_z ) = \mr{tr}(\lambda _2 ^* \mu \lambda _1)
	\]
	is  a Riemann form which defines a principal polarization $\rho _z$ on $A_z$. Also, 
	each element $a \in O$ gives a map $\theta(a) $ which maps $\CC ^2 \to \CC^2$
	and sends $\Lambda _z$ to itself, thereby inducing an endomorphism $\iota (a)$ of
	$A_z$. See \cite{Varieties-book}. We can enhance this structure by including a rigidification 
	of the points of order $N$ for some integer $N \ge 1$, namely an isomorphism
	\[
	\alpha :A[N] = \mr{Ker}(N : A \to A) \to (\ZZ/N)^4
	\]
	which carries Weil pairing to the standard symplectic pairing on $ (\ZZ/N)^4$. The resulting moduli 
	problem for $(A, \rho, \iota, \alpha)$ is now representable if $N \ge 3$, so we obtain fine moduli 
	schemes $X_B (N)$. Complex-analytically this is $\Gamma (N)\backslash \mathfrak{H}$ where
	$\Gamma (N)= \{\gamma \in \Gamma   \mid \gamma \equiv 1 \text{\ mod\ \ }N \}$.
	
	For the case of a general $T$, there are coarse moduli spaces $X_{\Gamma, T}$ that
	represent equivalence classes of $(A, \rho, \iota)$ for an equivalence relation arising 
	from stable quadratic twisting rings $R_T \subset O$. We refer to \cite{Rotger04}
	for details. For each prime $p$ dividing the discriminant $D$, there are involutions
	$w_p$ (Atkin-Lehner involutions) of the curve $X_{B}$. In fact
	$\Gamma _T = \langle \Gamma, w_p, p\in T \rangle$ and viewing the $w_p$ as 
	transformations of the moduli problem, 
	\[
	w_p (A, \rho, \iota) = (A/\mr{Ker}(I_p), \rho ', \iota ')  
	\]
	where $I_p\subset O$ is the set of elements whose norm is divisible by $p$. This 
	is a 2-sided ideal with $O/I_p \cong \FF_{p^2}$ and $I_p ^2 = p O$. Since
	$\mr{Ker}(I_p)$ is isotropic under the Weil pairing, $A/\mr{Ker}(I_p)$ inherits a principal 
	polarization $\rho '$, and $I_p$ being a 2-sided ideal, it also inherits an action $\iota '$
	of $O$. In other words, $X_{\Gamma, T} $ parametrizes $(A, \rho,  \iota)$ up to 
	identification of $A$ with $A/\mr{Ker}(I)$ where $I = \cap_{p \in T}  I_p$ for 
	$T \subset \Sigma$.
	
	Since these curves parametrize families of principally polarized abelian varieties 
	of dimension 2, there are modular embeddings
	\begin{align*}
		&f: \mathfrak{H} \to  \mathfrak{H}_2, \quad \mathfrak{H}_2 = \text{Siegel space of genus \ }2. \\
		& \varphi : O^1 \to \mr{Sp}_4(\ZZ), \quad f(\gamma z) =\varphi(\gamma )(f(z))
	\end{align*}
	inducing maps $X_{\Gamma, T} \to \mathfrak{A}_2$ to the moduli space of principally-polarized 
	abelian varieties of dimension 2 (see \cite{Hash95}). Because the Torelli map from the 
	moduli space of genus 2 curves
	$\mathfrak{M}_2 \to \mathfrak{A}_2 $ is a birational injection, we get algebraic coordinates 
	on  $\mathfrak{A}_2 $ from Igusa-Clebsch invariants, \cite{Igusa60}. The images of Shimura curves can then 
	be described in terms of Igusa-Clebsch invariants. See \cite{LinYang20} where Borcherds forms are used to calculate these.

	\subsection{ $D=6.$ The curves and modular forms}
	\label{SS:D=6.3}
	References for this section: \cite{Baba-Granath-g2}, \cite{Bayer-Travesa}.
	There is a unique quaternion algebra $\sf{B}_6$  of discriminant $D=6$ defined over $\QQ$. It has
	generators $1$, $I$, $J$, $K$ with $I^2=3$, $J^2=-1$, $K = IJ = -JI$. We fix the embedding 
	$\Phi : \sf{B}_6 \to M_2(\RR)$
	\[
	\gamma = x+yI + zJ + tK \mapsto \begin{bmatrix}x+y \sqrt{3} & z+t \sqrt{3}\\
		-( z-t \sqrt{3}) & x-y \sqrt{3}
	\end{bmatrix}.
	\]  
	All maximal orders are conjugate and we fix the representative 
	\[
	O_6 = \ZZ[1, I, J, (1+I+J+K)/2]. 
	\]
	The group $O_6^1$  can be 
	identified with its image $\Gamma _6 \subset \slr$ under $\Phi$.  
	Then the image of this in $\mr{PSL}_2(\RR)$ is $\bar{\Gamma}_6 = \Gamma _6 /\{\pm 1\}$
	and is a Fuchsian group without parabolic elements. 
	We let $X^{an}_6 =  \bar{\Gamma}_6 \backslash\mfr{H}$,  $X^{an}_6  = X_6(\CC)$ for a smooth projective curve defined over $\QQ$, Shimura's canonical model.  This is a 
	compact Riemann surface of genus 0. 
	In fact,  $X_6$ is the projective conic
	$x^2+3y^2 + z^2=0$.	
	
	The curve $X_6$ sits in a Galois $(\ZZ/2)^2$ tower of curves, defined by Atkin-Lehner quotients with $X^{(d)}_6:=X_6/\langle w_d\rangle$ and $X^{(+)}_6:=X_6/\langle w_2, w_3\rangle$. All these curves have genus zero and the corresponding groups have fundamental half-domains in the upper half-plane (or unit disk)
	which are hyperbolic polygons with either 3 or 4 sides. The angles at the vertices are $(\pi/a,\pi/b,\pi/c,(\pi/d))$ where $(a, b, c, (d))$ are as follows:

	$$
	\begin{diagram}
		&& X_6  & & \\
		& \ldTo ^ {\pi_2} & \dTo  _{\pi _3} & \rdTo  ^{\pi _6}\\
		X^{(2)}_6 & & X^{(3)}_6 & & X^{(6) }_6 \\
		& \rdTo & \dTo& \ldTo\\
		&& X ^{(+)}_6 &&
	\end{diagram}
	\qquad \qquad 
	\begin{diagram}
		&& \G^{(+)}_6:  (2,4,6)  & & \\
		& \ldLine  & \dLine   & \rdLine  \\
		\G^{(2)}_6: (3,4,4) &\quad & \G^{(3)}_6:(2,6,6) &\quad & \G^{(6) }_6:(2,2,2,3) \\
		& \rdLine & \dLine& \ldLine\\
		&& \G_6:(2,2,3,3) &&
	\end{diagram}
	$$

	\begin{figure}[h] \label{d6}
		\begin{center}
			\includegraphics[width=8.25cm,height=5.25cm]{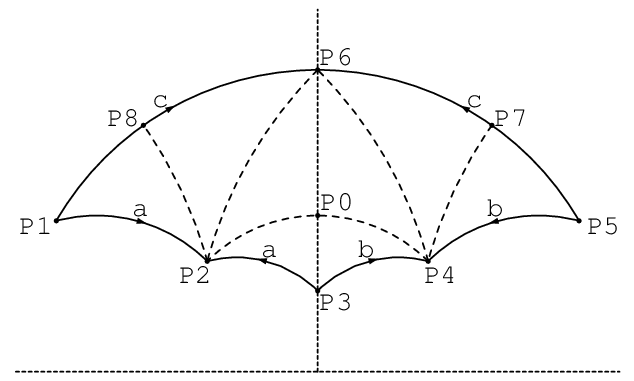}
			\caption{Fundamental domain for $\G_6$. Image from Bayer-Travesa, \cite{Bayer-Travesa}.  
				Atkin-Lehner quotients appear as dotted subregions. Special CM points and gluing of boundaries 
				shown.}
		\end{center}
	\end{figure}
	
	For the triangle groups listed above, we give a corresponding hypergeometric function. By suitable choices of two solutions to the corresponding differential equations (HDE), the ratios   give the
	Schwarz conformal maps from the upper half-plane to hyperbolic triangles:
	\vskip .2 cm
	\begin{center}
		\begin{tabular}{| c |c | c | }
			\hline
			curve &  angles & HDE \\
			\hline
			$X _6 ^{(2)}$ &  ( 3, 4, 4) & $ \pFq 21{\frac 1{12}&\frac1{3}}{&\frac{2}{3}}{t}$  \\
			\hline
			$X _6^{(3)}$ &   (2, 6, 6) & $ \pFq 21{\frac 1{12}&\frac1{4}}{&\frac{5}{6}}{t}$  \\
			\hline
			$X _6^{(+)}$ &  (2, 4, 6) & $ \pFq 21{\frac 1{24}&\frac7{24}}{&\frac{5}{6}}{t}$ \\
			\hline
		\end{tabular}
	\end{center}
	
	\begin{figure}[h]
		\begin{center}
			\includegraphics[width=5.25cm,height=5.25cm]{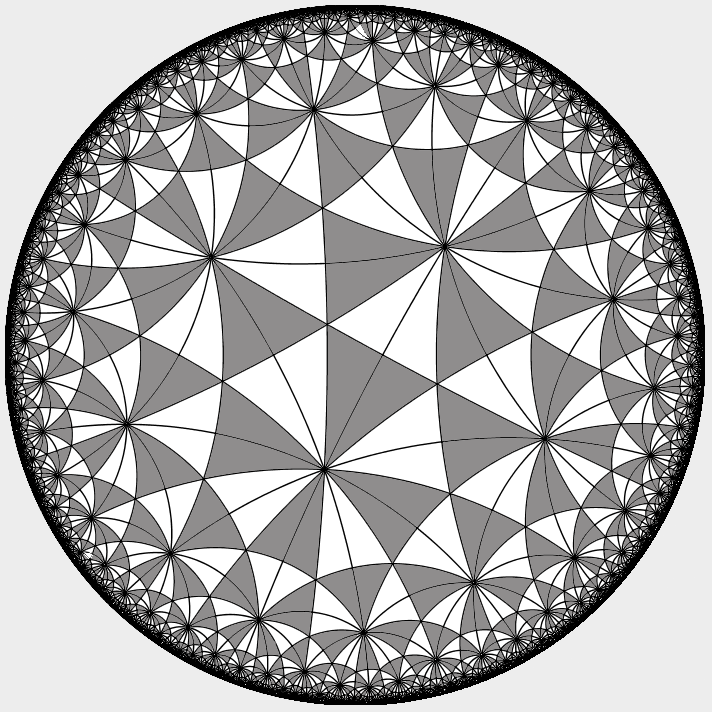}
			\caption{334 hyperbolic tessellation. Image thanks to  Helena Verrill.}
		\end{center}
	\end{figure}

	Let $\Gamma = O_6^1 $. The graded ring of modular forms is 
	\[
	\bigoplus _{k = 0} ^{\infty} S_{2k}(\Gamma)
	= \CC[ h_4(z), h_6(z), h_{12}(z)       ]
	\cong
	\CC[ h_4, h_6, h_{12}       ]/(  h_{12}^2 + 3 h_6 ^4 + h_4 ^6  )
	\]
	where $z $ is the coordinate in $\mathfrak{H}$, the subscript indicates the weight, and $h_4(z)$, 
	$h_6(z)$ are algebraically independent.	The Hauptmoduln on various 
	Atkin-Lehner quotients  are expressible in terms of these modular functions. For instance
	the $(2, 4, 6)$ curve we have denoted by $X_6 ^{+}$ is canonically isomorphic with 
	$\mathbb{P} ^1 _{\QQ}$ with a coordinate denoted $j_6$. 
	\[
	j_6 = \frac{16 h_6 ^4 }{9 h_4 ^6}.
	\]
	In the canonical model of the curve $X_6$ as the projective plane conic 
	$x^2 + 3y^2+z^2 = 0$ the Atkin-Lehner involutions are given on
	$p = (x, y, z)$ by $w_2 (p) = (x, -y, z)$,  $w_3 (p) = (-x, y, z)$, 
	$w_6 = w_2 w_3$. The curves $X_6 ^{(2)}$,  $X_6 ^{(3)}$, $X_6 ^{(6)}$, 
	$X_6 ^{+}$
	are all isomorphic with $\mathbb{P} ^1 _{\QQ}$, and the projections 
	$\pi _2$, $\pi _3$, $\pi _6$, $\pi _{+}$ are given by
	\[
	\pi _2 (p) = (x, z), \quad \pi _2 (p) = (y, z), \quad\pi _6 (p) = (x,  y), \quad
	\pi _{+}(p)  = (x^2, y^2)
	\]
	where  $(y/x)^2 = 9j_6/16 = h_6 ^4/h_4 ^6$.   See \cite{Baba-Granath-g2} for detailed discussion.
	\subsection{ $D = 6.$ Models of the abelian varieties}
	\label{SS:D=6.4}	
	If $\Gamma \subset \Gamma _6$ is a congruence subgroup with no torsion, then there is a universal 
	family of abelian varieties $f: A_{\Gamma}\to X_{\Gamma}$. In fact, the groups in our paper all have torsion
	so that we do not have universal families. Nonetheless we do have families of abelian varieties with quaternion structure
	that are ``sufficiently close'' to universal, and these are utilized for explicit calculations. This is analogous to the situation 
	for $\slz$, where no universal family of elliptic curves exists, but there are families of elliptic curves whose 
	$j$-invariant = $j$, the coordinate on the modular curve $X(\slz) = \mathbb{P} ^1 _j$: 
	\begin{equation}\label{eq:Universal-E}
		\mathcal E_j:\quad  y^2+xy=x^3-\frac{36x+1}{j-1728}.
	\end{equation} 
	These curves are not unique, but as Ihara and Scholl \cite{Sch88} have shown, they can be effectively used 
	for calculation of Hecke traces. There are three specific models of 2-dimensional abelian varieties with 
	multiplication by the quaternion algebra $B_6$ relevant to this paper:
	\begin{itemize}
		\item[1.] The Jacobian of the generalized Legendre curve 
		\[
		y^6 = x^4 (x-1)^3(x-\lambda)
		\]
		decomposes according to the characters of $\mu_6$. The part belonging to the primitive characters 
		is 2-dimensional. It corresponds to a hypergeometric motive and has QM by $B_6$. See \cite{WIN3a}.
		\item[2.] The subfamily of the Picard family of genus 3 curves
		\[
		y^3 = x (x-1)(x-\lambda)(x-\mu)
		\]
		with $\mu = 1-\lambda$ has the following property: Their Jacobians factor as $E \times A$ where 
		$E$ is an elliptic curve with CM by $\QQ (\zeta _3)$ and $A$ is a 2-dimensional abelian variety with endomorphisms by 
		$B_6$, \cite{Petkova-Shiga}

		\item[3.] The Baba-Granath curves. In  \cite{Baba-Granath-g2}, Baba and Granath construct a family of genus-2 curves whose Jacobians have QM by $B_6$. 
		This family is analogous the the family of elliptic curves with $j$-invariant = $ j$. 
		It is defined on a quadratic covering of $X(\Gamma _6) = \mathbb{P} ^1 _{j_6}$, not on  $\mathbb{P} ^1 _{j_6}$. 
		This is because of Mestre's obstruction: in general you cannot define a genus 2 curve $C$ over the same field 
		as the moduli of its Jacobian. Let $K = \QQ (j)$, $j = j_6$, and define $s = \sqrt{-6j}$, $t= -2(27j+16)$. Let
		$C$ be the projective nonsingular model of the curve $y^2 = f(x)$ where
		\begin{align*}
			f(x) = &(-4+3s)x^6 + 6tx^5 + 3t(28+9s)x^4 -4t^2 x^3\\
			&+3t^2(28-9s)x^2 + 6 t^3 x -t^3(4+3s). 
		\end{align*}
		We can rewrite this in terms of modular forms $h_4$, $h_6$, and $h_{12}$ in the previous section.
		Let $z\in \mathfrak{H}$. Then the curve $C_z: y^2 = g_z(x)$ where
		\begin{align*}
			g_z(x) =  &h_4^3 (z) (x^6 -21x^4 -21 x^2+1)\\
			&   +    \sqrt{-6}h_6 ^2 (z) (x^6 +9x^4 -9 x^2-1)  \\
			& + 2 \sqrt{2} h_{12}(z)x (3x^4 -2 x^2 +3) 
		\end{align*}
		has Jacobian isomorphic to $(A_z, \rho _z, \iota_z)$. The Igusa-Clebsch 
		invariant of this is
		\[
		[A, B, C, D] = [j_6+1, j_6, j_6(1-j_6), j_6^3]
		\]
		where 
		\[
		j_6 = \frac{AB -C}{AB+C} = \frac{D^2}{B^5}.
		\]

	\end{itemize}

	\bigskip

	\section{Main Theorem.}
	\label{S:main}
	
	\subsection{Outline.}
	\label{SS:main1}
	To prove the main claim,  formula \eqref{E:main}, by combining the trace formula \eqref{E:GL} with 
	the result of Eichler-Shimura theory \eqref{E:ES} we must calculate, for the groups we are 
	considering,
	\begin{equation*}
		\label{E:main2}
		-\mr{Tr}\, (T(p) \mid S_{2k}(U)) = 
		-\mr{Tr}\,   (\mr{Frob}_p \mid  H^1 _{et}(X_U \otimes _{\FF _p } \bar{\FF} _p,    V^{2k} (U)_{\ell})) 
		=  \sum _{x \in X_U(\FF_p)}\mr{Tr} (\mr{Frob}_x \mid        V^{2k} (U)_{\ell, \bar{x}}  ).
	\end{equation*}
	The contributions are of two sorts: (1) from cusps and elliptic points and (2) the rest. 
	\begin{itemize}
		\item[1.] The sheaf
		$V^{2k} (U)_{\ell}$ can be related recursively to $V^{2m} (U)_{\ell}$ for 
		$1 \le m < k$. This means that the traces  of Frobenius at $x$ on it can be expressed as polynomial 
		functions of the trace on $V^{2} (U)_{\ell, \bar{x}}$. It suffices therefore to relate the Frobenius trace 
		on $V^{2} (U)_{\ell, \bar{x}}$ to hypergeometric character sums for $x$ not a cusp or elliptic point. In fact one shows that
		\[
		V^{2} (U)_{\ell} \cong \text{a character of finite order}\otimes \text{a Kummer sheaf}
		\otimes \text{a hypergeometric sheaf\ }\mc{H}(HD (U))
		\]
		where $HD(U)$ is hypergeometric data attached to $U$.
		\item[2.] $V^{2} (U)_{\ell}$ (and $V^{2} (U)_{\CC}$)  is a rigid local system. Its local monodromies 
		are easily computed complex-analytically, and this is the same as the local monodromy $\ell$-adically. 
		The local monodromy of the hypergeometric sheaves are also known. The Kummer twist is there to insure 
		matching. Rigidity then shows that there is an isomorphism 
		\[
		V^{2} (U)_{\ell} \cong \text{a Kummer sheaf}
		\otimes \text{a hypergeometric sheaf\ }\mc{H}(HD (U))
		\]
		for the curves/sheaves over the algebraic closure $\bar{\QQ}$.
		
		\item[3.] To get an isomorphism over $\QQ$ we compare both sides at (one or more) CM point. For this 
		we need the explicit models of the curves and abelian varieties. 
		
		\item[4.] The Frobenius traces of the hypergeometric sheaves at points not corresponding to 
		cusps or elliptic points are given by hypergeometric character sums. For cusps or elliptic points, 
		this is essentially what is done in determining the fiber at places of bad reduction for a family 
		of abelian varieties. For elliptic curves, this is Tate's algorithm; for abelian varieties of dimension 2 coming
		from genus 2 curves, this is Qing Liu's algorithm \cite{QingLiu}. Actually, our situation is greatly simplified and we follow the 
		method of Scholl in \cite{Sch88}.

	\end{itemize}

	\subsection{Examples}
	\label{SS:main2}
	For the  triangle groups $\Gamma$ under consideration, we associate hypergeometric data $HD(\G)=\{\alpha(\G),\beta(\G)\}$ defined over $\QQ$
	for which Beukers-Cohen-Mellit \cite{BCM} introduced hypergeometric character sums $H_p(HD(\G), \lambda)$ for 
	$\lambda \in\FF_p^\times$. For each integer $m \ge 1$, let $F_m(S,T)$ be a degree-$m$ polynomial which expresses the symmetric polynomial 
	$ \sum_{i=0}^{2m} u^iv^{2m-i}$ in $u, v$ of degree $2m$ as a polynomial in $S=u^2+uv+v^2$ and $T=uv$, i.e., 
	\begin{equation}\label{eq:F}
		F_m(u^2+uv+v^2,uv)=
		\sum_{i=0}^{2m} u^iv^{2m-i}.
	\end{equation}
	
	\begin{theorem}[\cite{HLLT}]\label{thm:traceformula} For $\G = ~(2,4,6),~(2,\infty,\infty),(2,3,\infty),(2,4,\infty),(2,6,\infty)$, the table below describes the hypergeometric datum $HD(\G)=\{\alpha(\G),\beta(\G)\}$ and the choice of a generator $\l=\l(\G)$ of the field of $\QQ$-rational functions on $X_\G$ 
		by its values at each elliptic point of given order and cusp: 
		$$
		\begin{tabular}{|c|c|c|c|c|c|c|c|c|c|c|c|}
			\hline
			$\G$&$(2,\infty,\infty)$&$(2,3,\infty)$&$(2,4,\infty)$&$(2,6,\infty)$&$(2,4,6)$\\
			\hline
			$\l$&$(1,0,\infty)$&$(1,\infty, 0)$& $(1,\infty,0)$& $(1,\infty,0)$& $(-3,\infty,0)$\\ \hline
			$\alpha(\G)$    & $\{\f12,\f12,\f12\}$ &$\{\f12,\frac16,\frac56\}$&$\{\f12,\frac14,\frac34\}$&$\{\f12,\frac13,\frac23\}$&$\{\f12,\frac14,\frac34\}$\\
			\hline
			$\beta(\G)$&$\{1,1,1\}$&$\{1,1,1\}$&$\{1,1,1\}$&$\{1,1,1\}$&$\{1,\frac56,\frac76\}$\\ \hline
		\end{tabular}
		$$

		Given $k\ge 2$  an even integer and a fixed prime $\ell$, the terms on the right-hand side of (\ref{E:main2}) for almost all primes $p\ne \ell$ where $X_\G$ has good reduction are as follows. 
		
		For $\l$ not corresponding to an elliptic point or a cusp 
		\begin{equation}\label{eq:2}
			\mr{Tr}(\mr{Frob}_\l \mid (V^k(\Gamma)_{\ell})_{\bar{\lambda}} )=F_{k/2}(a_\G(\l,p),b_\G(\l,p)),
		\end{equation}
		where  $b_\G(\l,p)=p$  and
		\begin{equation}\label{eq:a_G}
			a_\G(\l,p)
			=\begin{cases}\left(\frac{1-1/\l}{p}\right)H_p(HD(\Gamma), 1/\lambda) & \text{ if $\G \ne (2,4,6)$};\\
				\left(\frac{-3(1+3/\l)}{p}\right)pH_p(HD(\G), -3/\lambda)& \text { if $\G=(2,4,6)$}.
			\end{cases}
		\end{equation}
		The contribution of $\l$ corresponding to a cusp
		is  $1$. There are also explicit formulas for the contributions of the elliptic points. 
	\end{theorem} 
	
	Here is an example illustrating our main result. For $\G=(2,4,6)$ the lowest $k$ with nontrivial $S_{k+2}(2,4,6)$ is  $k=6$, in which case $S_8(2,4,6)=\langle h_4^2 \rangle$ is $1$-dimensional, where $h_4$ generates $S_4(2,2,3,3)$. By Jacquet-Langlands correspondence \cite{JL}, $h_4^2$ corresponds to the normalized weight-$8$ level $6$ cuspidal newform $f_{6.8.a.a}$ in the LMFDB label.  For primes $p >5$, denote by $a_p(h_4^2)$ the eigenvalue of $T_p$ on $h_4^2$, which is equal to the $p$th Fourier coefficient $a_p(f_{6.8.a.a})$ of $f_{6.8.a.a}$. The theorem above gives
	
	\begin{align*}
		- a_p(h_4^2) = -a_p(f_{6.8.a.a})&=\sum_{\l \in \FF_p,\l \neq 0,-3} \left(a_{\G}(\l,p)^3-2pa_{\G}(\l,p)^2-p^2a_\G(\l,p)+p^3\right) \\    
		&+p((pH_p(HD(\G);1))^2-p^2) +\left(\left(\frac{-1}p\right)+\left(\frac{-3}p\right)+\left(\frac{-6}p\right)\right)p^3.
	\end{align*}The term $p((pH_p(HD(\G);1))^2-p^2)$ on the right side comes from the contribution at the elliptic point of order $2$ on $X_{(2,4,6)}$; it also equals $pa_p(f_{24.5.h.b})$, in terms of the $p$th Fourier coefficient of the weight-5 CM  modular form $f_{24.5.h.b}$.
	
	\subsection{Future directions}
	\label{SS:main3}	
	Here are two:
	\begin{itemize}
		\item[1.] Ohta's results apply to Shimura curves of quaternion algebras over any totally real number field $F$. 
		The $(5, 5, 5)$ triangle group,  which appears in Takeuchi's list, corresponds to a quaternion algebra $B$
		over $F = \QQ (\sqrt{5})$. The new part of the Jacobian of the curve
		\[
		y^{10} = x^2(1-x)^7 (1-\lambda x)^7
		\]
		is $4= 2[F:\QQ]$-dimensional and has endomorphism algebra the quaternion algebra over  $\QQ (\sqrt{5})$
		with discriminant $\mathfrak{p}_5$, the unique prime over $5$. This is a model of the 
		family of abelian varieties over a covering of the Shimura curve for $B$. See \cite{WIN3a}.
		The Jacquet-Langlands correspondence in this case will land in Hilbert modular forms 
		for  $\QQ (\sqrt{5})$.
		\item[2.] The quaternion algebra $B_{10}$ over $\QQ$ of discriminant 10 is also discussed in the 
		paper of Baba-Granath. But this time, the fundamental domains are not related to triangle groups. 
		In fact, the corresponding differential equation now has 4 regular singular points, and in particular 
		is no longer rigid. There is an accessory parameter which was determined by Elkies in \cite{Elkies}
		using Schwarzian equations. The number theory here will be governed by Heun functions over finite 
		fields, a theory that needs to be developed. See the entry on Heun Functions in \cite{DLMF}.

	\end{itemize}


	\bibliographystyle{plain}
	\bibliography{ref}

\end{document}